\documentclass[a4paper,11pt]{article}
\usepackage{pdfsync}
\usepackage{graphicx}
\usepackage{float}
\usepackage{color}
\usepackage{amsfonts,amscd,amssymb,amsmath,amsbsy,amsthm}
\usepackage{enumerate}
\usepackage{hyperref}
\usepackage[margin=1in]{geometry}
\usepackage{xspace}

\usepackage{lineno}

\newtheorem{theorem}{Theorem}[section]

\newtheorem{lemma}[theorem]{Lemma}

\newtheorem{mclaiminner}{Claim}
\newenvironment{mclaim}[1]{%
  \renewcommand\themclaiminner{#1}%
  \mclaiminner
}{\endmclaiminner}

\theoremstyle{definition}

\newtheorem{example}[theorem]{Example}

\theoremstyle{remark}

\numberwithin{equation}{section}

\newcommand*{\proofclaimtext}{Proof}
\newenvironment{proofclaim}[1][\proofclaimtext]{\begin{proof}[#1]}{\end{proof}}

\newenvironment{packed_enum}{
    \begin{enumerate}
        \setlength{\itemsep}{1pt}
        \setlength{\parskip}{0pt}
        \setlength{\parsep}{0pt}
}{\end{enumerate}}

\newenvironment{packed_itemize}{
    \begin{itemize}
        \setlength{\itemsep}{1pt}
        \setlength{\parskip}{0pt}
        \setlength{\parsep}{0pt}
}{\end{itemize}}

\newcommand{\heading}[1]{\medskip\par\noindent{\bf #1}}

\def\gS{\mathbb{S}} \def\gC{\mathbb{C}} \def\gA{\mathbb{A}} \def\gD{\mathbb{D}}

\def\calA{{\mathcal A}} \def\calB{{\mathcal B}} \def\calC{{\mathcal C}}
  \def\calF{{\mathcal F}}
\def\calG{{\mathcal G}}  
 \def\calK{{\mathcal K}} 
\def\calM{{\mathcal M}}  \def\calO{{\mathcal O}}
\def\calP{{\mathcal P}}

\def\planar{\hbox{\rm \sffamily PLANAR}\xspace}
\def\tree{\hbox{\rm \sffamily TREE}\xspace}

\newcommand{\pstab}[2]{{#1_{(#2)}}}
\newcommand{\sstab}[2]{{#1_{#2}}}
\newcommand{\gwr}{\mathrel{\wr\wr}}
\newcommand{\bo}{\partial}
\renewcommand{\int}{\mathring}

\DeclareMathOperator\Fix{Fix}
\DeclareMathOperator\Aut{Aut}
\DeclareMathOperator\Sym{Sym}
\DeclareMathOperator\Ker{Ker}
\DeclareMathOperator\id{id}
\DeclareMathOperator\Stab{Stab}

\title{Jordan-like characterization of automorphism groups of planar graphs\thanks{
We would like to thank Gareth Jones and Ilia Ponomarenko for several comments and suggestions that helped to improve this paper significantly.
The project was supported by GA\v{C}R 20-15576S. In addition, Roman Nedela was supported by APVV-15-0220, and Peter Zeman was supported by GAUK 124120.}}

\author{Pavel Klav\'{\i}k
	\and Roman Nedela
	\and Peter Zeman
}

\author{
  Pavel Klav\'{i}k
	\thanks{
  Department of Mathematics, University of West Bohemia, Pilsen, Czech Republic.
	\texttt{klavik@orgpad.com}.}
  \thanks{
  Orgpad, \url{https://www.orgpad.com}
  }
	\and
  Roman Nedela
	\thanks{
  Department of Mathematics, University of West Bohemia, Pilsen, Czech Republic.
	\texttt{nedela@savbb.sk}.}
	\and
  Peter Zeman
	\thanks{Department of Applied Mathematics,
	Faculty of Mathematics and Physics,
	Charles University, Prague, Czech Republic,
	\texttt{zeman@kam.mff.cuni.cz}.}
}
\date{}

\begin{document}

\maketitle

\begin{abstract}
We investigate automorphism groups of planar graphs.
The main result is a complete recursive description of all abstract groups that can be realized as automorphism groups of planar graphs.
The characterization is formulated in terms of inhomogeneous wreath products.
In the proof, we combine techniques from combinatorics, group theory, and geometry.
This significantly improves the Babai's description (1975).
\heading{Keywords:} planar graph, group, spherical group, wreath product.
\heading{MSC2020 classification:} primary 05C10, secondary 20D25, 05E18.
\end{abstract}

\section{Introduction}
\label{sec:intro}

The study of symmetries of discrete structures is one of the central topics in modern mathematics.
All such structures can be canonically encoded in by graphs~\cite{hedrlin_pultr}.
This motivates the study of symmetries of graphs.
Frucht (1939)~\cite{frucht} proved that every finite group can be realized as the automorphism group of some finite graph.
The following general problem arises.
For a family of graphs $\calC$, characterize the family $\Aut(\calC)$ of all abstract groups $G$ such that $G\cong \Aut(X)$, for some $X\in \calC$. We say that $\calC$ is \emph{universal} if $\Aut(\calC)$ contains all finite groups.
For a given class of graphs $\calC$, two important problems arise:
\begin{packed_itemize}
\item
to determine whether $\Aut(\calC)$ is universal,
\item
to describe which abstract groups belong to $\Aut(\calC)$.
\end{packed_itemize}

In this paper we investigate the automorphism groups of planar graphs.
The central problem considered in this paper reads as follows:
\medskip
\begin{quote}
\emph{Which abstract groups can be realized as the automorphism groups of planar graphs?}
\end{quote}
\medskip

The problem to characterize $\Aut(\calC)$ for distinguished graph classes $\calC$ have some history.
The family of all finite graphs is universal by the aforementioned result of Frucht~\cite{frucht}.
The bipartite and chordal graphs~\cite[Theorem 5]{lueker1979linear} are other well-known universal families of graphs.

On the other hand, the following very elegant description of the set of automorphism groups of trees as the smallest set of groups, closed under direct products and under wreath products with symmetric groups, is stated by Babai in~\cite[Proposition 1.15]{babai1996automorphism}.
Some of the ideas appear already in a paper by Jordan~\cite{jordan} published in 1869, therefore we call this and similar inductive characterizations \emph{Jordan-like}.
Denote the family of all finite trees by \tree. It follows that \tree is not a universal class of graphs.

Recently, $\Aut(\calC)$ was characterized for several families $\calC$ of graphs.
For instance, it was proved that interval graphs have the same automorphism groups as trees~\cite{kz,kz_full} (also follows from~\cite{algebraic_forests}), and that the unit interval graphs have the same automorphism groups as the disjoint unions of caterpillars~\cite{kz,kz_full}. For permutation and circle graphs~\cite{kz_full}, there are similar Jordan-like inductive descriptions as in Theorem~\ref{thm:jordan_trees}.
Comparability and function graphs form universal families even for the order dimension at most four~\cite{kz_full}.

Understanding the structure of $\Aut(\calC)$ may lead to an effective algorithm for computation of generators of $\Aut(X)$ for a graph $X \in \calC$.
The famous \emph{graph isomorphism problem} asks, whether two given graphs $X$ and $Y$ are \emph{isomorphic}.
This problem clearly belongs to NP\ and it is a prime candidate for an intermediate problem between P\ and NP-complete.
The graph isomorphism problem isrelated to computing generators of the automorphism groups.
Suppose that $X$ and $Y$ are connected graphs.
Then $X \cong Y$ if and only if some generator in $\Aut(X\mathbin{\dot\cup}Y)$ swaps $X$ and $Y$.
Thus, an algorithm for computing a generating set of the automorphism group of a graph gives an algorithm for solving the isomorphism problem.
On the other side, Mathon~\cite{mathon1979note} proved that generators of the automorphism group of an $n$-vertex graph can be computed using $O(n^3)$ instances of graph isomorphism.

\heading{Babai's Characterization.}
The automorphism groups of planar graphs were first described by Babai~\cite[Corollary
8.12]{babai1975automorphism} in 1973. The automorphism groups of 3-connected planar graphs can be
easily described geometrically by using the theorems of Whitney~\cite{whitney1932congruent} and
Mani~\cite{mani1971automorphismen}. The groups are called \emph{spherical groups} since they correspond to
the finite subgroups of the group of isometries of the sphere.
The automorphism groups of
$k$-connected planar graphs, where $k < 3$, are constructed by wreath products of the automorphism
groups of $(k+1)$-connected planar graphs and of stabilizers of some vertices and edges in the
automorphism groups of $k$-connected graphs.

In~\cite[p. 1457--1459]{babai1996automorphism}, Babai gives an overview of the key idea of a
reduction to the 3-connected case. Babai states there a consequence:

\begin{theorem}[Babai~\cite{babai1996automorphism}] \label{thm:short_babai}
If $X$ is planar, then the group $\Aut(X)$ has a subnormal chain
$$\Aut(X) = G_0 \triangleright G_1 \triangleright \cdots \triangleright G_m = \{1\}$$
such that each quotient group $G_{i-1} / G_i$ is either cyclic or symmetric or $\gA_5$.
\end{theorem}

This theorem explains what are the building blocks of the automorphism groups of planar graphs. To
understand how the group of automorphisms of a planar graph is composed from the building blocks is
a difficult problem. In particular, the extension of the action of a stabilizer of a vertex in the
automorphism group of a planar graph by a spherical group $\Sigma$ depends on the action of $\Sigma$
on vertices of the associated $3$-connected planar graph. To describe it in detail, one
needs to analyse the structure of point-orbits of distinguished spherical groups.
On the top of that Babai points out in~\cite[p.~69]{babai1975automorphism} that his characterization is not inductive:
``For the case of planar graphs, we determine the groups occurring in the Main Theorem, as abstract
groups (up to isomorphism). [\dots] It cannot be, however, considered as a characterization by
recursion of the automorphism groups of the planar graphs, since the group construction refers to
the action of the constituents of the wreath products.''

\heading{Our Characterization.}
The characterization of the automorphism groups of planar graphs, described in this paper, is inductive. The proof
uses three main ideas/concepts representing three fields of mathematics: graph theory, group theory and geometry.

The first key idea, is the idea of \emph{3-connected reduction}.
The reduction can be viewed as a function associating to a given planar graph $X$ an irreducible planar graph $X^R$ which is either $3$-connected, cycle, $K_2$ or $K_1$.
The crucial feature of the reduction is that information
on the automorphism group is preserved, i.e. the automorphism group of the original graph can be reconstructed from $X^R$ applying the reverse procedure, called expansion.
The idea of the 3-connected
reduction was first introduced in the seminal papers by Mac Lane~\cite{maclane} and
Trakhtenbrot~\cite{trakhtenbrot}, however. It was further extended
in~\cite{tutte_connectivity,quadratic_isomorphism_planar,
hopcroft_tarjan_dividing,cunnigham_edmonds,walsh,bienstock}.
The related decomposition (onto 3-connected components) can be represented by a tree whose nodes are 3-connected graphs, and this tree is
known in the literature mostly under the name \emph{SPQR-tree}~\cite{spqr1,spqr2,spqr3,spqr_linear}.

The main difference, compared to our approach, is that the reduction used in the former literature
applies exclusively to 2-connected graphs. We introduce a modified reduction
of~\cite{fkkn,fkkn16} which reduces simultaneously the parts separated by 1-cuts and 2-cuts. This
allows one to control the changes of the symmetries in each elementary reduction. In particular, the reduction process is done in a way that essential
information on the symmetries is preserved (using colored and directed edges), so that the
reconstruction of the automorphism group is possible.
A similar reduction was used in~\cite{fkkn,fkkn16} to study the behaviour of semiregular subgroups of $\Aut(X)$ with respect to 1-cuts and 2-cuts.

The reduction proceeds as follows. In each step, called an \emph{elementary reduction}, we replace all
atoms of the prescribed type of the considered graph $X$ by colored (possibly directed) edges, where the \emph{atoms} are
certain inclusion minimal subgraphs of $X$. This gives a reduction series consisting of graphs $X =
X_0,\dots,X_r$, where $X_{i+1}$ is created from $X_i$ by replacing all of its atoms with some edges.
The final graph $X_r$, called \emph{primitive}, contains no atoms. We show that
$X_r$ is either (essentially) 3-connected, or a cycle, or $K_2$, or $K_1$. We think of $X_r$ as a skeleton
associated to $X$: the graph $X$ is obtained from $X_r$ by expanding its edges to certain subgraphs
of $X$ separated by $1$- and $2$-cuts.

Another innovative feature is the way of coding the isomorphism classes of atoms during the reduction. Instead of using colors for both vertices and edges we extend the standard definition of a graph with multiple edges. This extension (by adding so called pendant edges) allow us to work with the vertex- and edge-coloring in a unified way.
For such general graphs we have to re-examine carefully the basic concepts such as the $n$-connectedness, for instance.
Here, Tutte's approach~\cite[Chapter 4]{tutte_graph_theory} turned out to be useful.

The next crucial concept is group-theoretical, the \emph{inhomogeneous wreath product}.
Similarly to the standard wreath product, it is a particular form of a semidirect product $K\rtimes H$, where the normal part $K$ factorises into a direct product of groups. Unlike the standard wreath product, these groups are not necessarily isomorphic. The outer group $H$ acts on the factors of $K$, permuting the factors, but it acts trivially on each of the factors. If all the factors of $K$ represent the same group $G$, then the inhomogeneous wreath product is just the standard wreath product $G\wr H$.

The inhomogeneous wreath products are introduced and investigated in Section~\ref{sec:groups}. Lemmas~\ref{lem:semidirect_by_direct},~\ref{lem:inhomogeneous_fixed},~\ref{lem:inhomogeneous_equivariant}
are of particular importance, since they allow us to
reduce a potentially unbounded number of group operations
to a few operations used in the main Theorems~\ref{thm:stabilizers} and~\ref{thm:composition}. It transpires that by employing the concept of
the inhomogeneous wreath product we are able to express the
main result in a comprehensive way (compare to~\cite[Pages 69-70]{babai1975automorphism}, see also Appendix~\ref{app:babai}).

The third fundamental tool is of geometric nature.
Our classification relies on determining the \emph{equivariance classes of point orbits} in
the actions of spherical groups done in Section~\ref{sec:spherical}. In the analysis we employ the
fruitful concept of a quotient orbifold introduced by Thurston in his investigation of 3-dimensional manifolds.

The above mentioned
ideas are employed in Section~\ref{sec:characterization}, where the main results of the paper are proved. First, in Theorem~\ref{thm:stabilizers} , we present
an inductive characterization of the vertex-stabilizers of planar graphs. The basic groups are cyclic, dihedral and
symmetric groups. The set of abstract groups that appear as the vertex-stabilizers of vertices of planar graphs is characterized as the smallest set of groups containing the basic groups and closed with respect to five types of inhomogenous wreath products. Theorem~\ref{thm:composition} completes the characterization, the resulting set of groups is obtained
from the set of vertex-stabilisers using the inhomogeneous
wreath products with the spherical groups. In particular,
each spherical group determines one inhomogeneous wreath product. It is interesting that our characterization describes the automorphism groups of planar
graphs without referring to planarity explicitly, as a simple recursive process which builds them
from the basic groups.

\section{Graphs}
\label{sec:graphs}

Although our primary objective is to investigate automorphism groups of simple planar graphs, our reduction procedure naturally produces parallel edges and edges that are incident to only one vertex.
Therefore, for technical reasons, we adopt an extended model of a graph inspired by topological graph theory, where graphs are understood as
1-dimensional CW-complexes.

Usually, a 1-dimensional CW-complex
is formed by a set of disjoint  $1$-simplices homeomorphic
to open intervals, called edges, and a discrete set of points,
called vertices. In the standard model, the closure of each edge
is homeomorphic either to a closed interval bounded by two vertices, or to a simple circle containing exactly one vertex (such edges are called loops).
In our model, we further allow that some of the vertices are missing.
Moreover, for technical reasons we exclude loops.
We describe a combinatorial model of such graphs.

\heading{Extended graphs.}
An \emph{extended multigraph} $X$ (or just a \emph{graph} $X$) is a $3$-tuple $(E, V, \iota)$, where $E$ is a finite set of edges, $V$ is a finite set of vertices with $V\cap E=\emptyset$, and $\iota : E \to \binom{V}{2} \cup \binom{V}{1} \cup \{\emptyset\}$ is an incidence function.
An edge $e \in E$ is \emph{standard} if $|\iota(e)| = 2$, \emph{pendant} if $|\iota(e)| = 1$, and \emph{free} if $|\iota(e)| = 0$.
A pendant edge $e$ with $\iota(e) = \{v\}$ is called a \emph{ single pendant edge} if for every other pendant edge $e'$ with $\iota(e') = \{u\}$ we have $v\neq u$.

To avoid confusion when working with more graphs, we use $V(X)$ and $E(X)$ to denote the vertices and edges of $X$, respectively.
For any $v \in V$, we define the \emph{degree} of $v$, denoted by $\deg(v)$, to be the number of edges incident to $v$. 

We denote the $n$-cycle, for $n\geq 3$, by $C_n$.
A \emph{dipole} $D_n$ is a subgraph consisting of two vertices joined by $n\geq 2$ parallel edges.
A \emph{pendant star} $S_n$ is a subgraph having one-vertex incident to $n\geq 2$ pendant edges.
Every pendant edge is either a single pendant edge, or it is contained in a pendant star.

\heading{Subgraphs of extended graphs.}
A graph $X' = (E', V', \iota')$ is a subgraph of $X = (E, V, \iota)$ if $E' \subseteq E$, $V' \subseteq V$, and $\iota'(e) = \iota(e) \cap V'$ for every $e \in E'$.
For every $v \in V$, we define the \emph{degree} of $v$ in $X'$, denoted by $\deg_{X'}(v)$, to be the number of edges of $X'$ incident to $v$.

For $S \subseteq V(X)  \cup E(X)$, we denote by $X \setminus S$ the subgraph with the set of edges $E(X)\setminus S$ and with the set of vertices $V(X)\setminus S$.
In particular, note that if $S \subseteq V(X)$, then $X\setminus S$ contains the edges incident to the vertices in $S$.

We say that a set $U$ of vertices of $X$ is an \emph{$n$-cut} if there are subgraphs of $Y$ and $Z$ of $X$ such that $X = Y\cup Z$, $V(Y)\cap V(Z) = U$, $|U| = n$, $E(Y)\cap E(Z) = \emptyset$, and both $Y$ and $Z$ have at least $n$ edges.
A $1$-cut is often called an \emph{articulation}.
A graph $X$ is \emph{$m$-connected} if it has no $n$-cut, for every $n<k$.
There are several definitions of $m$-connectedness, this one is due to Tutte~\cite[Page 70]{tutte_graph_theory}.
A maximal $2$-connected subgraph of a graph $X$ is called a \emph{block}.

For a subgraph $A$ of $X$, we use the topological notation to denote the \emph{boundary}
$\bo A$  and the \emph{interior} $\int A$ of $A$.  We set $\bo A$ equal to the set of vertices of
$A$ which are incident with an edge not contained in $A$. The interior of $A$ is defined by setting $\int A = A \setminus \bo A$.

\heading{Colored graphs with directed edges.}
We often consider graphs with colored edges. Moreover, some of the edges of the considered graphs may  be directed.
The reductions replace parts of the graph by new edges and the colors encode
isomorphism classes of replaced parts.
Even if the input graph $X$ is simple, colored multigraphs are naturally constructed in the reduction process.
Note that in the model of edge-colored multigraphs we use, one can realize a vertex-coloring by attaching colored single pendant edges to vertices.

\heading{Homomorphisms of extended graphs.}
Let $X=(E, V, \iota)$ and $X'=(E', V', \iota')$ be two graphs.
A \emph{homomorphism} $X\to X'$ is a mapping $\alpha : E\cup V \to E'\cup V'$ such that $\alpha(E)\subseteq E'$, $\alpha(V)\subseteq V'$, and
$$
\alpha(\iota(e)) = \iota'(\alpha(e)) \quad \text{for every } e \in E.
$$

For colored graphs with directed edges, a homomorphism preserves both the colors and directions of edges.
An \emph{isomorphism} $X$ and $X'$ is a bijective homomorphism $X \to X'$.
If $X$ and $X'$ are isomorphic, we write $X\cong X'$.
An \emph{automorphism} of $X$ is an isomorphism of $X$ to itself.
We denote the \emph{automorphism group} of $X$ by $\Aut(X)$. 

For subgraphs of $X$, we usually consider only isomorphisms preserving their
boundaries. Let $A$, $A'$ be subgraphs of $X$. An isomorphism $\alpha \colon A \to A'$ is called a
\emph{$\bo$-isomorphism} if $\alpha(\bo A) = \bo \alpha(A)$. If such a $\bo$-isomorphism exists, we say that
$A$ is $\bo$-isomorphic to $A'$, denoted $A \cong_\bo A'$.
Observe that for every subgraph $A$ and every automorphism $\alpha$ of $X$, the mapping $\alpha |_A$ is a $\bo$-isomorphism from $A$ to $\alpha(A)$.

\section{Elements of group theory}
\label{sec:groups}

In this section, we introduce the notation and concepts from group theory, which we use throughout the paper.
Most importantly, we introduce the inhomogeneous wreath product, which is a quite general tool allowing us to describe the automorphism group of a graph in terms of the automorphism groups of its smaller subgraphs.
We conclude the section with two preliminary applications of the inhomogeneous wreath product.

In the paper, we use the following notation for some standard families of groups:
\begin{packed_itemize}
\item $\gS_n$ is the \emph{symmetric group} of all permutations of the set $\{1,\dots,n\}$,
\item $\gC_n$ is the \emph{cyclic group} of integers $\{0,\dots,n-1\}$ with addition modulo $n$,
\item $\gD_n$ is the \emph{dihedral group} of the symmetries of a regular $n$-gon, and
\item $\gA_n$ is the \emph{alternating group} of all even permutations of the set $\{1,\dots,n\}$.
\end{packed_itemize}
We note that $\gD_1 \cong \gC_2$, $\gD_2 \cong \gC_2^2 = \gC_2\times \gC_2$, $\gC_3\cong \gA_3$, and $\gD_3 \cong \gS_3$. 

\subsection{Group actions}

A group $G$ \emph{acts} on a set $\Omega$ if there is a mapping $\cdot : G \times \Omega
\to \Omega$ such that for every $x\in\Omega$ and for every $g,h\in G$, we have
$1 \cdot x = x$ and $(gh) \cdot x = g \cdot (h \cdot x)$.
If there is no confusion, we write $gx$ instead of $g\cdot x$.
The set $\Omega$ is called a \emph{$G$-set} if $G$ acts on $\Omega$.
For $x \in \Omega$, we define the \emph{stabilizer} $G_x = \{g : g\cdot x = x\}$.
For $S \subseteq \Omega$, we define the \emph{set-wise stabilizer} $G_S = \{g : g\cdot S = S\}$ and the \emph{point-wise} stabilizer $G_{(S)} = \{g : gx = x, x \in S\}$.

For $x \in \Omega$ we define the \emph{orbit} $[x]_G = \{gx : g \in G\}$.
We write $[x]$ if $G$ is apparent from the context.
The action of $G$ on $\Omega$ is \emph{transitive} if it acts on $\Omega$ with a single orbit, \emph{semiregular} if $G_x$ is trivial for every $x\in\Omega$, and \emph{regular} if it is transitive and semiregular.

\heading{Equivariant actions.}
Two $G$-sets $\Omega_1$ and $\Omega_2$ are \emph{equivariant} if there is a bijection $f : \Omega_1 \to \Omega_2$ such that $f(g\cdot x) = g\cdot f(x)$ for every $g\in G$ and $x \in \Omega_1$.
Equivariance is an equivalence relation on the set of all $G$-spaces and we call its classes \emph{equivariance classes}.
If $\Omega$ is a $G$-set, then we say that two orbits $[x]_G$ and $[y]_G$ are \emph{equivariant orbits} if the induced $G$-sets $[x]$ and $[y]$ are equivariant.

\begin{example}
\label{ex:two_actions_dihedral}
We define two non-equivariant actions of $\gD_n$, to which we refer repeatedly in the paper.
Let $\Omega = \{0,\dots,n-1\}$ be the \emph{vertices} of a regular $n$-gon, for $n\geq 2$, and let $\Omega' = \{\{i,i+1 \mod n\} : i = 0,\dots,n-1\}$ be the \emph{edges} of the regular $n$-gon.
We define \emph{the action of $\gD_n = \langle r,t \mid r^n = (rt)^2 = 1\rangle$ on the vertices of the $n$-gon} by setting $r\cdot i = i + 1 \mod n$, for $i=0,\dots,n-1$, and $t\cdot i = n-i$, for $i = 0,\dots,n - 1$.
We define \emph{the action of $\gD_n$ on the edges of the $n$-gon} to be the action on $\Omega'$ induced by the action of $\gD_n$ on $\Omega$.
If $n$ is even, then, by Lemma~\ref{lem:dih_actions}, these two actions are non-equivariant, although they are isomorphic.
\qed
\end{example}

The following lemma is well-known; see for example~\cite{cameron}.

\begin{lemma}
\label{lem:equivariant_action_characterization}
Let $\Omega$ be a $G$-set, and let $G/G_x$, for $x\in\Omega$, be the set of all the left cosets of $G_x$.
Then $[x]_G$ and $G/G_x$ are equivariant $G$-sets.
Moreover, two $G$-sets $[x]_G$ and $[y]_G$ are equivariant if and only if $G_x$ and $G_y$ are conjugate in $G$.
\end{lemma}

\begin{lemma}
\label{lem:trivial_orbits_equivariant}
If $G$ acts on $\Omega$, then all orbits of size $|G|$  are equivariant and all orbits of size $1$ are equivariant.
\end{lemma}

\begin{proof}
Let $x,y \in \Omega$.
If $[x]$ and $[y]$ are of size $|G|$, then both  $G_x$ and $G_y$ are trivial.
If $[x]$ and $[y]$ are of size $1$, then both $G_x$ and $G_y$ equal to $G$.
In both cases the statement immediately follows from Lemma~\ref{lem:equivariant_action_characterization}.
\end{proof}

\begin{lemma}
\label{lem:dih_actions}
Let $G \cong \gD_n = \langle r, t \mid r^n = t^2 = (rt)^2 = 1\rangle$ act on a set $\Omega$.
If $n$ is odd, then all orbits of size $n$ are equivariant.
If $n$ is even, then there are three equivariance classes of orbits of size $n$ determined by the subgroups $\langle r^{n/2} \rangle$, $\langle t\rangle$, $\langle rt\rangle$, respectively.
\end{lemma}

\begin{proof}
By the orbit-stabilizer theorem, $G_x$ of point is isomorphic to $\gC_2$.
If $n$ is even, there are exactly three conjugacy classes of subgroups of order $2$.
If $n$ is odd there is just one conjugacy class of $G_x$ in $G$.
The statement follows from Lemma~\ref{lem:equivariant_action_characterization}.
\end{proof}

\subsection{Group products}

Given two groups $K$ and $H$, and a group homomorphism $\theta \colon H \to \Aut(K)$ taking $h\mapsto \theta_h$, the \emph{outer semidirect product} $K \rtimes_{\theta} H$ is the cartesian product of $K$ and $H$ with  the operation is defined by the rule
$$(g_1, h_1) \cdot (g_2,h_2) = (g_1\theta_{h_1}(g_2), h_1h_2).$$
Note that if $\theta$ is trivial, then $K\rtimes_\theta H = K\times H$ is the \emph{(outer) direct product}.
A group $G$ is an \emph{inner semidirect product} if $G = KH$, where $K\lhd G$, $H\leq G$, and $K\cap H = \{1\}$.

It is well-known~\cite{rotman} that inner and outer semidirect products are equivalent in the following sense. For every inner semidirect product $G=KH$, there is a homomorphism $\theta\colon H \to \Aut(K)$, defined by $\theta_h(x) = hxh^{-1}$, such that $G\cong K\rtimes_\theta H$.
On the other hand, if $G = K\rtimes_\theta H$ is an outer semidirect product, then $G = K'H'$, where $K' = \{(k,1) : k \in K\}$ and $H' = \{(1,h) : h \in H\}$.
Thus, we use both inner and outer semidirect product interchangeably.

In this paper, we use only particular special cases of semidirect products, which we now introduce.

\heading{Wreath product.}
Let $D$ and $Q$ be groups, let $\Omega$ be a finite $Q$-set, and let $K = \prod_{\omega\in\Omega}D_\omega$, where $D_\omega\cong D$ for all $w\in \Omega$.
Then the \emph{wreath product} of $D$ by $Q$, denoted by $D\wr Q$, is the semidirect product of $K$ by $Q$, where $Q$ acts on $K$ by $q\cdot (d_w) = (d_{q\omega})$ for $q \in Q$ and $(d_w) \in \prod_{w\in\Omega}D_\omega$.
The normal subgroup $K$ of $D\wr Q$ is called the \emph{base} of the wreath product.

\heading{Inhomogeneous wreath product.}
We generalize the standard wreath product as follows.
Let $D_1,\dots,D_m$ and $Q$ be groups, let $\Omega$ be a finite $Q$-set with orbits $\Omega_1,\dots,\Omega_m$, and let $K = \prod_{i=1}^m\prod_{\omega \in \Omega_i}D_{i,\omega}$, where $D_{i,\omega} \cong D_i$ for $\omega \in \Omega_i$ and $i = 1,\dots,m$.
Then the \emph{inhomogeneous wreath product} of the groups $D_1,\dots,D_m$ by $Q$, denoted by $(D_1,\dots,D_m)\gwr Q$, is the semidirect product of $K$ by $Q$, where $Q$ acts on $K$ by $q\cdot (d_{i,w}) = (d_{i,q\omega})$ for $q \in Q$ and $(d_{i,w}) \in \prod_{i=1}^m\prod_{w\in\Omega_i}D_{i,\omega}$.
The normal subgroup $K$ of $(D_1,\dots, D_m) \gwr Q$ is called the \emph{base} of the inhomogeneous wreath product.
Note that by choosing $D_1\cong \cdots\cong D_m \cong D$ we obtain the wreath product $D\wr Q$.

If $\Lambda_i$ is a $D_i$-set, for $i=1,\dots,m$, then $\bigcup_i\Lambda_i\times\Omega_i$ can be made into a $((D_1,\dots,D_m)\gwr Q)$-set.
Given $d \in D_i$ and $\omega \in \Omega_i$, for $i = 1,\dots,m$, we define the permutation $d_{i,\omega}^*$ of $\bigcup_i\Lambda_i\times\Omega_i$ as follows:
for each $(\lambda, \omega^\prime) \in \bigcup_i\Lambda_i\times\Omega_i$, set
$$
d_{i,\omega}^*(\lambda,\omega^\prime) =
\begin{cases}
(d\lambda, \omega^\prime) \quad &\text{if} \quad \omega^\prime = \omega,\\
(\lambda, \omega^\prime) &\text{if} \quad \omega^\prime \neq \omega.
\end{cases}
$$
It is easy to see that $d_{i,w}^*d_{i,w}^{\prime *} = (dd^\prime)_{i,\omega}^*$.
Thus,
$$D_{i,w}^* = \{d_{i,w}^* : d \in D_i\}$$
is a subgroup of $\Sym(\bigcup_i\Lambda_i\times\Omega_i)$.
This is certified by the map $D_i \to D_{i,\omega}^*$, given by $d \mapsto d_{i,\omega}^*$, which is an isomorphism.
For each $q \in Q$, we define a permutation $q^*$ of $\bigcup_i\Lambda_i\times\Omega_i$ by
$$q^*(\lambda,\omega^\prime) = (\lambda, q\omega^\prime),$$
and define
$$Q^* =\{q^* : q \in Q\}.$$
It is easy to see that $Q^*$ is a subgroup of $\Sym(\bigcup_i\Lambda_i\times\Omega_i)$.
This is certified by the map $Q \to Q^*$, given by $q \mapsto q^*$, which is an isomorphism.

We defined the inhomogeneous wreath product as an outer semidirect product.
We note that an analogical concept for permutation groups was defined in~\cite{babai1975automorphism,robinson} under the names ``$G$-Kranz product'' and ``generalized composition'', respectively.
The following theorem, in which we characterize it as an inner semidirect product, is a generalization of~\cite[Theorem 7.24]{rotman}.

\begin{theorem}
\label{thm:inhomogeneous_wr}
Let $D_1,\dots,D_m$ and $Q$ be groups, let $\Omega$ be a $Q$-set with orbits $\Omega_1,\dots,\Omega_m$, and let $\Lambda_i$ be a $D_i$-set for $i = 1,\dots,m$.
Then the inhomogeneous wreath product $(D_1,\dots,D_m) \gwr Q$ is isomorphic to the subgroup
$$W = \langle \bigcup_{i=1}^{m}\bigcup_{\omega\in\Omega_i}D_{i,\omega}^*, Q^* \rangle \leq \Sym(\bigcup_{i=1}^m \Lambda_i\times\Omega_i).$$
In particular, $W = K^*Q^*$, where $K^* = \langle\bigcup_{i=1}^{m}\bigcup_{\omega\in\Omega_i}D_{i,\omega}^*\rangle$ is the direct product $\prod_{i=1}^m\prod_{\omega\in\Omega_i}D_{i,\omega}^*$.
\end{theorem}

\begin{proof}
We show that the group $K^*$ is the direct product.
It is easy to se that $d_{j,\omega^\prime}^{\prime *}d_{i,\omega}^*d_{j,\omega^\prime}^{\prime *-1} = d_{i,\omega}^*$ for all $d_{i,\omega}^*$ and $d_{j,\omega^\prime}^{\prime *} \in K^*$, i.e, $D_{i,\omega}^* \lhd K^*$ for all $i = 1,\dots,m$ and $\omega \in \Omega$.
Further, each $d_{i,\omega}^* \in D_{i,\omega}^*$ fixes all $(\lambda, \omega^\prime) \in \bigcup_j\Lambda_j\times\Omega_j$ with $\omega^\prime \neq \omega$, while each element of $\langle\bigcup_{j=1}^m\bigcup_{\omega^\prime\in\Omega_i,\omega^\prime\neq\omega}D_{j,\omega^\prime}^*\rangle$ fixes all $(\lambda, \omega)$ for all $\lambda \in \Lambda_i$.
It follows that if $d_{i,\omega}^* \in D_{i,\omega}^*\cap\langle\bigcup_{j=1}^m\bigcup_{\omega^\prime\in\Omega_i,\omega^\prime\neq\omega}D_{j,\omega^\prime}^*\rangle$, then $d_{i,\omega}^* = 1$.

If $d_{i,\omega}^* \in D_{i,\omega}^*$ and $q^* \in Q^*$, then $q^*d_{i,\omega}^*q^{*-1} = d_{i,q\omega}^*$.
Hence, $q^*K^*q^{*-1} \leq K^*$ for each $q^* \in Q^*$.
Since $W = \langle K^*, Q^*\rangle$, it follows that $K^* \lhd W$.
We get that $W = K^*Q^*$ and to see that $W$ is a semidirect product of $K^*$ by $Q^*$, it suffices to show that $K^*\cap Q^* = 1$.
If $d_{i,\omega}^* \in K^*$, then either $d_{i,\omega}^*(\lambda, \omega^\prime) = (d\lambda, \omega^\prime)$ or $d_{i,\omega}^*(\lambda, \omega^\prime) = (\lambda, \omega^\prime)$, i.e., each $d_{i,\omega}^*$ always fixes the second coordinate.
If $q^*\in Q^*$, then $q^*(\lambda,\omega^\prime) = (\lambda, q\omega^\prime)$, i.e., each $q^*$ always fixes the first coordinate.
Therefore every element of $K^*\cap Q^*$ fixes every $(\lambda, \omega^\prime)$ and hence must be equal to $1$.

Finally, the map $(D_1,\dots,D_m)\gwr Q \to W$, defined by $(d_{i,\omega})q \mapsto (d_{i,\omega}^*)q^*$, is an isomorphism.
\end{proof}

\begin{lemma}
\label{lem:semidirect_by_direct}
Let $G = KQ$ be a semidirect product of $K$ by $Q$.
Further, let $K = K_1\times K_2$ and $Q = Q_1\times Q_2$ such that $Q_1$ fixes $K_2$ pointwise and $Q_2$ fixes $K_1$ pointwise.
Then,
$$G \cong K_1 \rtimes Q_1 \times K_2 \rtimes Q_2.$$
\end{lemma}

\begin{proof}
First, we prove that $K_1Q_1 \lhd KQ$.
For $k_1q_1 \in K_1Q_1$ and $kq\in KQ$, we have
$$kqk_1q_1(kq)^{-1} = kqk_1q^{-1}qq_1q^{-1}k^{-1}.$$
There are $k_1^\prime \in K_1$ and $q_1^\prime \in Q_1$ such that the right-hand side is equal to
$$kk_1^\prime q_1^\prime k^{-1} = kk_1^\prime q_1^\prime k^{-1}q_1^{\prime -1}q_1^\prime.$$
Now, if $k \in K_1$, then $q_1^\prime k^{-1}q_1^{\prime -1} \in K_1$.
On the other hand if $k \in K_2$, then $q_1^\prime k^{-1}q_1^{\prime -1} = k^{-1}$ and the right-hand side is equal to $k_1^\prime q_1^\prime$.
Similarly, one can prove that $K_2Q_2\lhd KQ$.
It is easy to see that $K_1Q_1 \cap K_2Q_2 = \{1\}$.
Finally, $K_i\lhd Q_i$ and so $K_iQ_i$ is a semidirect product, for $i = 1, 2$.
\end{proof}

\begin{lemma}
\label{lem:inhomogeneous_fixed}
Let $D_1,\dots,D_m$ and $Q$ be groups, let $\Omega$ be a $Q$-set with orbits $\Omega_1,\dots,\Omega_m$.
Moreover, assume that all $\Omega_1,\dots,\Omega_k$ are of size one for some $k\leq m$.
Then,
$$(D_1,\dots,D_m) \gwr Q \cong D_1\times\cdots\times D_k \times (D_{k+1},\dots,D_m) \gwr Q.$$
\end{lemma}

\begin{proof}
By the definition, there are groups $D_{i,\omega}\cong D_i$, for $i=1,\dots,m$, such that $(D_1,\dots,D_m)\gwr Q$ is the semidirect product of $K$ by $Q$, where $K = \prod_{i=1}^m\prod_{\omega\in\Omega_i}D_{i,\omega}$ and $Q$ acts on $K$ by $q\cdot (d_{i,\omega}) = (d_{i,q\omega})$.
By the assumptions, we can write
$$K = D_{1,\omega_1}\times\cdots\times D_{k,\omega_k}\times\prod_{i=k+1}^m\prod_{\omega\in\Omega_i}D_{i,\omega},$$
where $\Omega_{i} = \{\omega_i\}$, for $i\leq k$.
By the definition of inhomogeneous wreath product, all subgroups $D_{i,\omega_i}$ are centralized by $Q$, which completes the proof.
\end{proof}

\begin{lemma}
\label{lem:inhomogeneous_equivariant}
Let $D_1,\dots,D_m$ and $Q$ be groups, let $\Omega$ be a $Q$-set with orbits $\Omega_1,\dots,\Omega_m$.
Moreover, assume that the orbits $\Omega_1$ and $\Omega_2$ are equivariant.
Then,
$$(D_1,\dots,D_m) \gwr Q \cong (D_1\times D_2, D_3,\dots,D_m) \gwr Q.$$
\end{lemma}

\begin{proof}
By the definition, there are groups $D_{i,\omega}\cong D_i$, for $i=1,\dots,m$, such that $(D_1,\dots,D_m)\gwr Q$ is the semidirect product of $K$ by $Q$, where $K = \prod_{i=1}^m\prod_{\omega\in\Omega_i}D_{i,\omega}$ and $Q$ acts on $K$ by $q\cdot (d_{i,\omega}) = (d_{i,q\omega})$.
Further, let $f\colon \Omega_1\to\Omega_2$ be a bijection certifying that $\Omega_1$ and $\Omega_2$ are equivariant orbits, i.e., $f(q\omega) = qf(\omega)$ for every $q\in Q$ and $\omega \in \Omega_1$.
We can write
$$K = \prod_{\omega\in \Omega_1}D_{1,\omega}\times D_{2,f(\omega)} \times \prod_{i=3}^m\prod_{\omega\in\Omega_i}D_{i,\omega}.$$
Let $d_{1,\omega}d_{2,\omega}^{\prime} \in D_{1,\omega}\times D_{2,f(\omega)}$ for some $\omega \in \Omega_1$.
From the definition of inhomogeneous wreath product and equivariance of the orbits $\Omega_1$ and $\Omega_2$ it follows that
$$qd_{1,\omega}d_{2,\omega}^\prime q^{-1} = qd_{1,\omega}q^{-1}qd_{2,\omega}^\prime q^{-1} = d_{1,q\omega}d_{2,qf(\omega)}^\prime = d_{1,\omega}d_{2,f(q\omega)}^\prime \in D_{1,q\omega}\times D_{2,f(q\omega)}.$$
Thus, the inhomogeneous wreath product $(D_1,\dots,D_m)\gwr Q$ is determined by the induced action of $Q$ on $\Omega\setminus \Omega_2$ and therefore it is isomorphic to $(D_1\times D_2, D_3,\dots,D_m)\gwr Q$.
\end{proof}

\subsection{Two simple applications}

We show how can the inhomogenous wreath product be used to describe the automorphism groups of several well-known classes of graphs.

\heading{Automorphism groups of disconnected graphs.}
There is a well-known description of the automorphism group of a disconnected graph in terms of the automorphism groups of its connected components.

\begin{theorem}
\label{thm:aut_disconnected}
Let $X_1, \dots, X_n$ be pairwise non-isomorphic simple connected graphs and let $X$ be the disjoint union $X = \bigcup_{i=1}^n\bigcup_{j=1}^{m_i}X_{i,j}$, where $X_{i,j}\cong X_i$. Then
$$\Aut(X) \cong \Aut(X_1) \wr \gS_{m_1} \times \cdots \times \Aut(X_n) \wr \gS_{m_n}.$$
\end{theorem}

\begin{proof}
First, we describe the automorphisms of $X$.
We denote $V(X_{i,j}) = \{(v, j) : v \in V(X_i), j \in \{1\dots,m_i\}\}$.

For every $d \in \Aut(X_i)$, there is an automorphism $d_{i,j}^* \in \Aut(X)$ such that $d_{i,j}^*(v,j^\prime) = (dv, j^\prime)$ if $j^\prime = j$, and $d_{i,j}^*(v,j^\prime) = (v, j^\prime)$ if $j^\prime \neq j$.
Note that the subgroup $D_{i,j}^* = \{d_{i,j}^* : d \in \Aut(X_i)\} \leq \Aut(X)$ fixes all the vertices in $V(X)\setminus V(X_{i,j})$ pointwise.

Further, for every $q \in \gS_{m_1}\times\cdots\times\gS_{m_n}$, there is an automorphism $q^* \in \Aut(X)$ of $X$ such that $q^*(v,j) = (v, qj)$.
Clearly the subgroup $Q^* = \{q^* : q \in \gS_{m_1}\times\cdots\times\gS_{m_n}\} \leq \Aut(X)$ acts on the set $\bigcup_i\{(i,j) : j = 1,\dots,m_i\}$ of the labels of the connected components of $X$.

We have $\Aut(X) = \langle \bigcup_i\bigcup_jD_{i,j}^*, Q^*\rangle$ and by Theorem~\ref{thm:inhomogeneous_wr}
$$\Aut(X) \cong (\Aut(X_1),\dots,\Aut(X_n))\gwr  (\gS_{m_1}\times\cdots\times\gS_{m_n}).$$
Finally, applying Lemma~\ref{lem:semidirect_by_direct} repeatedly on the right-hand side gives the theorem.
\end{proof}

\heading{Automorphism groups of trees.}
By a simple application of the previous theorem, we prove the Jordan's theorem~\cite{jordan} characterizing the automorphism groups of trees.

\begin{theorem}\label{thm:jordan_trees}
The class $\Aut(\tree)$ is equal to the class of groups $\calG$ defined inductively as follows:
\begin{packed_itemize}
\item[(a)]
$\{1\} \in \calG$.
\item[(b)]
If $G_1, G_2 \in \calG$, then $G_1 \times G_2 \in \calG$.
\item[(c)]
If $G \in \calG$, then $G \wr \gS_n \in \calG$ for all $n \in \mathbb{N}$.
\end{packed_itemize}
\end{theorem} 

\begin{proof}
First, we show that $\Aut(\tree)\subseteq\calG$.
We proceed by induction on the number of vertices.
Every tree has a \emph{center}, which is either a vertex or an edge.
For a tree $T$, the group $\Aut(T)$ stabilizes its center setwise.
Moreover, if the center of $T$ is an edge, we form $T'$ by subdividing this edge by exactly one vertex.
Then $\Aut(T')\cong \Aut(T)$ and the center of $T'$ is a vertex.
Thus, we may assume that the center is a vertex $v$.
After deleting $v$, we get a disconnected graph consisting of $k_i$ copies of the tree $T_i$, for $i = 1,\dots,n$.
By Theorem~\ref{thm:aut_disconnected}, we have
$$\Aut(T)\cong \Aut(T_1)\wr \gS_{m_1}\times\cdots\times\Aut(T_n)\wr\gS_{m_n},$$
and so, by the induction hypothesis, $\Aut(T)$ is isomorphic to a group in $\calG$.

To prove $\Aut(\tree)\supseteq\calG$, we proceed by the induction on the number of operations.
For each of the two operations, we easily construct a tree realizing the corresponding operation.
\end{proof}

\section{Reduction to 3-connected graphs} \label{sec:reduction}

We develop a reduction into 3-connected components.
Under certain conditions, we show that the automorphism group $\Aut(X)$ can be inductively reconstructed from automorphism groups of colored 3-connected components.
We stress that planarity of $X$ is not assumed in this section.

\subsection{Atoms}

\begin{figure}[t!]
\centering
\includegraphics{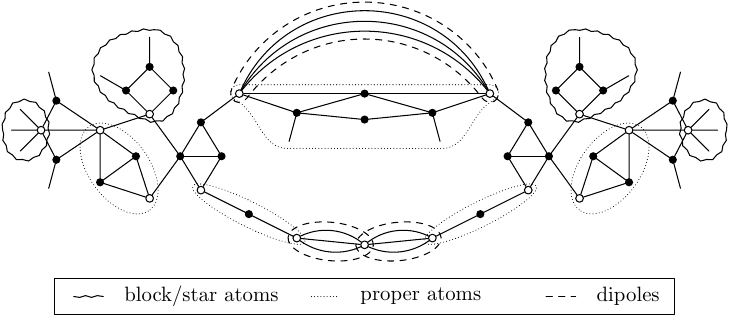}
\caption{An example of a graph with denoted atoms. The white vertices belong to the boundary of some
atom, possibly several of them.}
\label{fig:examples_of_atoms}
\end{figure}

We introduce the concept of atoms which are the basic components used in the reduction.
To this end, we say that a \emph{non-trivial $2$-cut of $X$} is $2$-cut $U = \{u,v\}$ of a block $B$ of $X$ such that $\deg_B(u) \geq 3$ and $\deg_B(v) \geq 3$.
We start by defining several types of subgraphs of a graph $X$, called \emph{parts}.
\begin{packed_itemize}
\item
A \emph{star part} is any maximal pendant star, and a \emph{dipole part} is any maximal dipole.
\item
Let $X = A\cup B$ be a decomposition defining a $1$-cut such that both $A$ and $B$ are not single pendant edges and $A$ is connected.
Then $A$ is a \emph{block part}.
\item
A subgraph $A$ is a \emph{proper part} if $X = A\cup B$ is a decomposition defining a non-trivial $2$-cut $U = \{u,v\}$ such that both $A$ and $A \setminus \{u,v\}$ are connected.
\end{packed_itemize}
Denote by $\calP(X)$ is the union of all parts of the graph $X$.
The inclusion-wise minimal elements of $\calP(X)$ are called \emph{atoms} and the set of all atoms of $X$ will be denoted by $\calA(X)$.
We distinguish \emph{star atoms}, \emph{dipole atoms}, \emph{block atoms}, and \emph{proper atoms}, according to the type of the defining part.
Observe that all star parts and dipole parts are already atoms; see Figure~\ref{fig:examples_of_atoms}. 

We conclude this section by stating several lemmas, which were proven in~\cite{3connected_regular}.
For completeness, we include their proof in Appendix~\ref{app:lemmas}.

\heading{Structure of primitive graphs and atoms.}
A graph containing no atoms is called \emph{primitive}.
The key idea of the reduction process will be to replace atoms by edges until we obtain a primitive graph.
The graph $\overline{X}$ obtained from $X$ by removing all the pendant edges is called the \emph{essence} of $X$.
We say that $X$ is \emph{essentially $3$-connected} if $X$ contains no pendant stars and  $\overline{X}$ is $3$-connected.
Similarly, $X$ is \emph{essectially a cycle} if $X$ has no pendant stars and $\overline{X}$ is a cycle.
Let $A$ be a proper atom with $\bo A = \{u,v\}$.
Note that, by the definition, in a proper atom, the vertices in $\bo A$ are always non-adjacent.
The \emph{extended proper atom $A^+$} is $A$ together with the edge $uv$.

\begin{lemma}
\label{lem:essence}
If $X$ is a primitive graph, or a block atom, or an extended proper atom, on least $3$ vertices, then the essence $\overline{X}$ of $X$ is either $3$-connected, or a cycle.
\end{lemma}

\begin{lemma}
\label{lem:structure_primitive}
If $X$ is a primitive graph, then $X$ is a $3$-connected graph, or $C_n$ for $n\geq 3$, or $K_2$, or $K_1$, or it can be obtained from these graphs by attaching single pendant edges. 
\end{lemma}

\begin{lemma}
\label{lem:structure_block}
Every block atom $A$ is $K_2$ possibly with a unique single pendant edge attached, or essentially a cycle, or essentially $3$-connected.
\end{lemma}

\begin{lemma}
\label{lem:structure_proper}
For every proper atom $A$, the extended atom $A^+$ is either essentially a cycle, or essentially $3$-connected.
\end{lemma}

The following lemma states that atoms can overlap only in their boundaries.
This is important for the reduction process to be uniquely determined in every step.

\begin{lemma}
\label{lem:interiors}
Let $A_1$ and $A_2$ be two distinct atoms in a connected graph. Then $A_1\cap A_2=\partial A_1\cap \partial A_2$, in particular $\int A_1\cap \int A_2=\emptyset$.
\end{lemma}

\heading{Action of the automorphism group on atoms.}
We observe that atoms behave well with respect to the action of the automorphism group, making it possible to replace them by edges, while preserving some information about the automorphism group of the graph.

\begin{lemma}
\label{lem:atoms_and_automorphisms}
Let $A$ be an atom and let $g \in \Aut(X)$.
Then the following hold:
\begin{packed_itemize}
\item
We have $g(A) \cong A$, $g(\bo A) = \bo g(A)$, and $g(\int A) = \int {g(A)}$.
\item
If $g(A) \neq A$, then $g(\int A) \cap \int A = \emptyset$ and $g(A)\cap A = \bo g(A)\cap \bo A$.
\end{packed_itemize}
\end{lemma}

Based on the previous lemma, an automorphism can either fix the boundary of a proper or dipole atom pointwise, or it can swap it.
This motivates the following definition.
An atom $A$ is

\begin{packed_itemize}
\item \emph{symmetric} if $\pstab{\Aut(A)}{\bo A}$ is a proper subgroup of $\sstab{\Aut(A)}{\bo A}$,
\item \emph{asymmetric}, if $\pstab{\Aut(A)}{\bo A} = \sstab{\Aut(A)}{\bo A}$.
\end{packed_itemize}
Note that star and block atoms are always symmetric.

\begin{figure}[b!]
\centering
\includegraphics{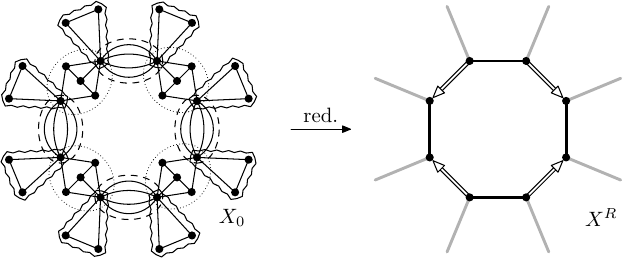}
\caption{On the left, the graph $X_0$ has three isomorphism classes of atoms, one of each type.  We
reduce $X_0$ to $X^{R}$ which is an eight cycle with single pendant edges, with four black symmetric
edges replacing the dipoles, four gray undirected edges replacing the block atoms, and four white
directed edges replacing the proper atoms.  The reduction series ends with $X^R$ since it is
primitive.}
\label{fig:example_of_reduction}
\end{figure}

\subsection{Reduction}

For a connected graph $X$, the reduction produces a series of graphs $X = X_0,\dots,X_r$.
To construct the graph $X_{i+1}$ from $X_i$, we first find the collection of all atoms $\calA = \calA(X_i)$ of $X_i$.
To each isomorphism class of $\calA$, we assign one new color not yet used in the graphs $X_0,\dots,X_i$.

The interior of every star atom and block atom with $\bo A = \{u\}$ will be replaced by a pendant edge $e_A$ of the assigned color.
The interior of every proper atom or a dipole atom $A$ with $\bo A = \{u,v\}$ by a new edge $e_A$ joining $u$ and $v$ of the assigned color.
If $A$ is symmetric, the edge $e_A$ is undirected. If $A$ is asymmetric, we chose a direction of $e_A$ and we extend the orientation consistently with the action of the automorphism group  to the entire isomorphism class of $e_A$.

Let $\calA_1(X_i)$ denote the set of all star, dipole and proper atoms and let $\calA_2(X_i)$ denote the set of all block atoms.
The reduction from $X_i$ to $X_{i+1}$ is done by applying one of the following two steps:
\begin{packed_enum}
\item
If $\calA_1(X_i) \neq \emptyset$, then we replace the interior of every atoms $A \in \calA_1(X_i)$.
\item
If $\calA_1(X_i) = \emptyset$, then we replace the interior of every atom $A \in \calA_2(X_i)$.
\end{packed_enum}

\begin{figure}[t!]
\centering
\includegraphics{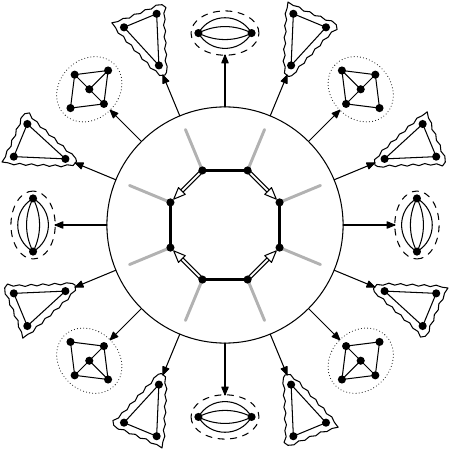}
\caption{An example of a reduction step.}
\label{fig:reduction_tree}
\end{figure}

According to Lemma~\ref{lem:interiors}, the replaced interiors of the atoms are pairwise disjoint, so the reduction is well-defined.
By applying the reductions repeatedly we obtain after $r$ steps a primitive graph $X_r$.
For a given graph $X$, we denote by $X^R = X_r$ the result of applying the reduction process to $X$.
For an example of a reduction step see Figure~\ref{fig:example_of_reduction} and Figure~\ref{fig:reduction_tree}.

\heading{Properties of the reduction.}
For a graph $X$, the incidence graph of blocks and articulations of the essence $\overline{X}$ of $X$ is called the \emph{block tree} of $X$.
The center of the block tree either corresponds to an articulation or to a block.
In the former case we say that $X$ has a \emph{central articulation} and in the latter case we say that $X$ has a \emph{central block}.

\begin{lemma}
Let $B$ be the central block or the central articulation of $X$.
Then $V(X^R) \subseteq V(B)$.
\end{lemma}

\begin{proof}
We show  that during each elementary reduction, the central block/articulation does not change in the block trees.
This is clearly true if $\calA_1(X_i) \neq \emptyset$ since in this case the block tree of $X_i$ and of $X_{i+1}$ coincide.
If $\calA_1(X_i) = \emptyset$, then the block tree of $X_{i+1}$ is constructed from that of $X_i$ by removing all the leaves and the vertices adjacent to them.
This does not change the position of the center.
\end{proof}


\heading{Expansions of edges.}
Consider the reductions series $X = X_0,\dots,X_r$.
For an edge $e \in E(X_i)$, for $i > 0$, let $A_e$ be either the atom of $X_{i-1}$ which reduces to $e$, or $A_e = e$ if $e \in E(X_{i-1})$.
The \emph{expansion} $A_e^{*}$ of $e \in E(X_i)$ is the unique subgraph of $X = X_0$ which reduces to $e$ after $i$ steps.
Observe that $\iota(e) = \bo A_e = \bo A_e^*$.

\heading{Reduction epimorphism.}
We investigate how the groups $\Aut(X_i)$ and $\Aut(X_{i+1})$ are related.
We define a mapping $\varphi_i\colon \Aut(X_i) \to \Aut(X_{i+1})$ as follows.
For $g \in \Aut(X_i)$, if $g(A_e) = A_f$, then we set $\varphi_i(g)(e) = f$.

\begin{lemma}
The mapping $\varphi_i\colon \Aut(X_i)\to\Aut(X_{i+1})$ is a group epimophism.
\end{lemma}

\begin{proof}
For $g,h \in \Aut(X_i)$, we have
$$\varphi_i(gh)(e_A) = e_{gh(A)} = e_{g(h(A))} = \varphi_i(g)(e_{h(A)}) = \varphi_i(g)\varphi_i(h)(e_A).$$
Hence, $\varphi_i$ is a group homomorphism.

We show that $\varphi_i$ is surjective.
For $g' \in \Aut(X_{i+1})$, we want to extend $g'$ to $g \in \Aut(X_i)$ such that $\varphi_i(g) = g'$.
We just describe this extension on a single edge $e \in E(X_{i+1})$.
If $e$ is an original edge of $X_i$, there is nothing to extend.
Suppose that $e = e_A$, for some atom $A$ of $X_i$.
Then $f = g'(e)$ is an edge of the same color and the same type as $e$, and therefore $f$ arises from $A_f \cong A_e = A$.
The automorphism $g'$ prescribes the action of $g$ on the boundary $\bo A$.
We need to show that it is possible to extend the action of $g'$ to $\int A$ consistently.
The proof splits depending on symmetry type of $A$.
\begin{packed_itemize}
\item
\emph{$A$ is a star atom or a block atom.}
Then edges $e$ and $f = g'(e)$ are pendant, incident to articulations $u$ and $v = g'(u)$.
We define $g_{\restriction \int A_e}$ to be a $\bo$-isomorphism $A_e \to A_f$.
\item
\emph{$A$ is an asymmetric proper atom or an asymmetric dipole.}
By the definition, the orientations of $e$ and $f = g'(e)$ is consistent with respect to $g'$.
Since $\int A_e \cong \int A_f$, we define $g_{\restriction \int A_e}$ to be an orientation-preserving $\bo$-isomorphism $A_e \to A_f$.
\item
\emph{$A$ is a symmetric proper atom or a symmetric dipole.}
Let $h\colon A_e\to A_f$ be a $\bo$-isomorphism.
If $h$ maps $\bo A_e$ consistently with $g'$, then we set $g_{\restriction \int A_e} = h_{\restriction \int A_e}$.
Otherwise, we set $g_{\restriction \int A_e} = h_{\restriction \int A_e}\circ k$, where $k$ is an automorphism of $A_e$ swapping the two vertices in $\bo A_e$.
\end{packed_itemize}
It follows that $\varphi_i$ is an epimorphism.
\end{proof}

By the First isomorphism theorem,
we know that $\Aut(X_i)$ is an extension of $\Ker(\varphi_i)$ by $\Aut(X_{i+1})$ and
$$\Aut(X_{i+1}) \cong \Aut(X_i) / \Ker(\varphi_i).$$
The following lemma gives a structural information on the kernel of $\varphi_i$.

\begin{lemma}
\label{lem:kernel_of_epimorphism}
The kernel $\Ker(\varphi_i)$ is isomorphic to the direct product $\prod_{e \in E(X_{i+1})} \pstab{\Aut(A_e)}{\bo A_e}$.
\end{lemma}
\begin{proof}
For an atom $A \in \calA(X_i)$, denote by $K_A$ the point-wise stabilizer of $X\setminus \int A$ in $\Ker(\varphi_i) \leq \Aut(X_i)$.
Clearly, $K_A \cong \pstab{\Aut(A)}{\bo A}$.
According to Lemma~\ref{lem:interiors}, the interiors of any two disctinct atoms $A,B \in \calA(X_i)$ are pairwise disjoint.
It follows that $K_A\cap K_B = \id$ and that both $K_A$ and $K_B$ are normal subgroups.
Hence, $\Ker(\varphi_i) = \prod K_A \cong \prod\pstab{\Aut(A)}{\bo A}$.
\end{proof}

We say that an atom $A$ with $\bo A = \{u,v\}$ is \emph{centrally symmetric} if there exists an automorphism $h\in \sstab{\Aut(A)}{\bo A}$ such that $h(u) = v$, $h(v) = u$, $h^2 = \id$, and $h$ centralizes $\pstab{\Aut(A)}{\bo A}$.
Note that all star and block atoms are centrally symmetric.
Also, every symmetric dipole is centrally symmetric.
Moreover, if each symmetric proper atom is centrally symmetric, then the following theorem hold.

\begin{theorem} \label{thm:semidirect_product}
Let $X_0,\dots,X_r$ be the reduction series of a graph $X$.
If every symmetric proper atom of $X_i$, is centrally symmetric, for $i < r$, then
$$\Aut(X_i) \cong \left (\prod_{e\in E(X_{i+1})} \pstab{\Aut(A_e)}{\bo A_e}\right )\rtimes \Aut(X_{i+1})$$
is the inhomogeneous wreath product defined by the action of $\Aut(X_{i+1})$ on $E(X_{i+1})$.
\end{theorem}

\begin{proof}
For simplicity, we denote $G_i = \Aut(X_i)$ and $K_i = \Ker(\varphi_i)$.
To prove the theorem, we first find a subgroup $H$ of $G_i$ such that $G_i = K_iH$, $K_i\cap H = \{\id\}$, and $H \cong \Aut(X_{i+1})$.

The idea of the proof is as follows.
We form a sequence of graphs $X_{i+1} = Y_0,\dots,Y_s = X_i$, where $s$ is the number of all edge-orbits of $G_{i+1}$ and  $Y_j$ is obtained from $Y_{j-1}$ by expanding all the edges of one edge-orbit into atoms.
First, we set $H_0 = G_{i+1}$.
For $j > 0$, let $\calO = [o]_{G_{i+1}}$ be the $j$th edge-orbit of $G_{i+1}$, where $o\in E(X_{i+1})$.
For $g \in G_{i+1}$, we find an extension $g' \in \Aut(Y_j)$ such that $g'(A_e) = A_f$ if and only if $g(e) = f$ for all $e, f \in\calO$.
The set $H_j = \{g' : g \in G_{i+1}\}$ is then a subgroup of $\Aut(Y_j)$ isomorphic to $G_{i+1}$.
It is easy to see that $H_j \cap K_i = \{\id\}$, for $j = 0,1,\dots,s$.
Finally, we put $H = H_s$.

We describe the outlined construction for every fixed $j > 0$.
Assume that we already constructed the groups $H_0,\dots,H_{j-1}$.
Let $\calO = [o]_{G_{i+1}}$ be the $j$th edge-orbit of $G_{i+1}$, as above.
To construct $H_j$, we distinguish several cases according to the type of $A_o$.

\begin{figure}[b!]
\centering
\includegraphics{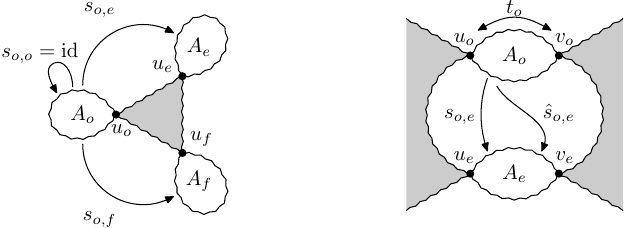}
\caption{Case 1 is demonstrated on the left, the respective block
atoms are $A_o$, $A_e$ and $A_f$. Case 3 is demonstrated on the right, the additional involution $t_e \in \sstab{\Aut(A_e)}{\bo A_e}$ transposes $u_e$
and $v_e$.}
\label{fig:group_extension}
\end{figure}

\begin{itemize}
\item
\emph{Case 1: The atom $A_o$ is a star or a block atom.}
All the edges in $\calO$  are pendant.
For $e \in \calO$, let $u_e$ be the articulation such that $\bo A_e = \{u_e\}$.
Choose arbitrarily $\bo$-isomorphisms $s_{o,e}$ from $A_o$ to $A_e$, and put $s_{o,o} = \id$ and $s_{e,f} = s_{o,f}s^{-1}_{o,e}$.
If $g'(e) = f$, we set $g_{\restriction_{\int{A_e}}} = {s_{e,f}}_{\restriction_{\int{A_e}}}$.
Since
\begin{equation} \label{eq:sigma}
s_{e,c} = s_{f,c} s_{e,f},\qquad \forall e,f,c \in \calO,
\end{equation}
the composition of the extensions $g_1$ and $g_2$ of two automorphisms $g'_1$ and $g'_2$ is defined on the interiors of all $\{A_e : e \in \calO\}$ exactly as the extension of $g'_2 g'_1$.
Therefore the extensions form a group.

\item
\emph{Case 2: The atom $A_o$ is an asymmetric proper atom or dipole.}
All the edges in the orbit $\calO$ are oriented consistently with the action of $G_{i+1}$, and the end-vertices form two orbits.
For $e \in \calO$, let $\iota(e) = \{u_e,v_e\}$, where $\{ u_e : e \in \calO\}$ and $\{v_e : e \in \calO \}$ are the two vertex-orbits.
Further, the construction proceeds as in Case 1, in addition, since $A_o$ is asymmetric, we have that $s_{o,e}(u_o) = u_e$ and $s_{o,e}(v_o) = v_e$.

\item
\emph{Case 3: The atom $A_o$ is a symmetric proper atom or a dipole.}
All the edges in the orbit $\calO$ are standard edges, and their end-vertices form one orbit. 
For $e \in \calO$, let $\iota(e_i) = \{u_e,v_e\}$, where $\{u_e, v_e : e \in \calO\}$ is the vertex-orbit.
Further, for $e \in \calO$ we arbitrarily choose one $\bo$-isomorphism $s_{o,e}$ from $A_o$ to $A_e$ such that $s_{o,e}(u_o) = u_e$ and $s_{o,e}(v_o) = v_e$, and set $s_{e,f} = s_{o,f} s^{-1}_{o,e}$, for $e,f\in\calO$.

By the assumptions, there is an involution $t_o \in \sstab{\Aut(A_o)}{\bo A_o}$ which exchanges $u_o$ and $v_o$.
Then $t_o$ defines an involution of $A_e$ by conjugation $t_e = s_{o,e} t_o s^{-1}_{1,e}$.
It follows that
\begin{equation}
\label{eq:tau}
t_f = s_{e,f} t_e s^{-1}_{e,f},\qquad\text{and consequently}\qquad s_{e,f} t_e = t_f s_{e,f},\qquad \forall e,f\in\calO.
\end{equation}
We put $\hat s_{e,f} = s_{e,f} t_e = t_f s_{e,f}$ which is an isomorphism mapping $A_e$ to $A_f$ such that $\hat s_{e,f}(u_e) = v_f$ and $\hat s_{e,f}(v_e) = u_f$.
In the extension, we put $g_{\restriction_{\int{A_e}}} = {s_{e,f}}{\restriction_{\int{A_e}}}$ if $g'(u_e) = u_f$, and $g_{\restriction_{\int{A_e}}} = \hat {s_{e,f}}_{\restriction_{\int{A_e}}}$ if $g'(u_e) = v_f$. 

Besides (\ref{eq:sigma}), we get the following additional identities:
\begin{equation} \label{eq:hatsigma}
\hat s_{e,c} = s_{f,c} \hat s_{e,f},\qquad
\hat s_{e,c} = \hat s_{f,c} s_{e,f},\quad\text{and}\quad
s_{e,c} = \hat s_{f,c} \hat s_{e,f},\qquad \forall e,f,c \in \calO.
\end{equation}
We just argue the last identity:
$$\hat s_{f,c} \hat s_{e,f} = t_c (s_{f,c} s_{e,f}) t_e = t_c
s_{e,c} t_e = t_c t_c s_{e,c} = s_{e,c}.$$
 It follows that the composition of the extensions
$g_2g_1$ is correctly defined, and it coincides with the extension of
$g'_2g'_1$.
\end{itemize}

From the construction of $H$ it follows that $H\cong \Aut(X_{i+1})$ and $H\cap K_i = \{\id\}$.
Therefore, $G_i = K_iH$ is a semidirect product.
Note that the complement $H$ of $\Ker(\varphi_i)$, in the statement, is not uniquely determined. 
In the proof, the group $H$ depends on the choice of $s_{1,i}$ and $t_1$.

We prove that the semidirect product $K_iH$ can be expressed as an inhomogeneous wreath product.
Let $A_e$ be a centrally symmetric atom with $\bo A_e = \{u_e,v_e\}$.
Recall that symmetric dipoles are also centrally symmetric.
In the case when each proper atom is centrally symmetric, we can choose $t_1 = t_e \in \sstab{\Aut(A_e)}{\bo A_e}$ to be the central involution.
Let $g \in \pstab{\Aut(A_e)}{\bo A_e}$.
We need to prove that for every $h \in H$ mapping $A_e$ to $A_f$, we have that $hgh^{-1} \in \pstab{\Aut(A_f)}{\bo A_f}$ does not depend on the choice of $h$.
Let $s_{e,f}$ be the $\bo$-isomorphism mapping $A_e \to A_f$, chosen above.
Now, if $s_{e,f}(u_e) = h(u_e)$, then $hg h^{-1} = s_{e,f}gs_{e,f}^{-1}$.
If $s_{e,f}(u_e) \neq h(u_e)$, then
$$hgh^{-1} = t_fs_{e,f}gs_{e,f}^{-1}t_f = s_{e,f}t_egt_es_{e,f}^{-1} = s_{e,f}gs_{e,f}^{-1}.$$
The second equality follows from~(\ref{eq:tau}) and the last equality holds since $t_e$ is a central involution.
\end{proof}

\subsection{Recursive construction of automorphism groups} \label{sec:reduction_aut_groups}

We show that the reduction can be used inductively to describe automorphism groups of graphs in terms of the automorphism groups of their 3-connected components.
Let $\theta_i = \varphi_{i-1} \circ \cdots \circ \varphi_0$ denote the epimorphism $\Aut(X_0)\to \Aut(X_i)$, for $i=1,\dots,r$.

\begin{lemma}\label{lem:kernel_product}
The kernel $\Ker(\theta_i)$ is isomorphic to the direct product $\prod_{e\in E(X_{i})} \pstab{\Aut(A_e^*)}{\bo A_e^*}$.
\end{lemma}
\begin{proof}
The statement follows by repeatedly applying Lemma~\ref{lem:kernel_of_epimorphism}.
\end{proof}

\begin{theorem}
\label{thm:expanded_inhomogeneous_product}
Let $X_0,\dots,X_r$ be the reduction series of a graph $X$.
If every symmetric proper atom of $X_0,\dots,X_i$ is centrally symmetric, for $i < r$, then
$$\Aut(X) \cong \left (\prod_{e\in E(X_{i+1})} \pstab{\Aut(A_e^*)}{\bo A_e^*}\right )\rtimes \Aut(X_{i+1})$$
is the inhomogeneous wreath product defined by the action of $\Aut(X_{i+1})$ on $E(X_{i+1})$.
\end{theorem}

\begin{proof}
By Theorem~\ref{thm:semidirect_product} the statement holds for $i=0$.
If $i>0$, we shall use Theorem~\ref{thm:semidirect_product} repeatedly, to conclude that $\Aut(X)$ contains an isomorphic copy of $\Aut(X_{i+1})$ forming a complement of $\Ker(\theta_{i+1})$.
\end{proof}

\begin{theorem}
\label{lem:semidirect_expatom}
Let $X_0,\dots,X_r$ be the reduction series of a graph $X$.
If every symmetric proper atom of $X_0,\dots,X_i$ is centrally symmetric, for $i < r$, and $f \in E(X_{i+1})$, then
$$\pstab{\Aut(A_f^*)}{\bo A_f^*} \cong \left (\prod_{e\in E(A_f)} \pstab{\Aut(A_e^*)}{\bo A_e^*} \right )\rtimes \pstab{\Aut(A_f)}{\bo A_f}$$
is the inhomogeneous wreath product defined by the action of $\pstab{\Aut(A_f)}{\bo A_f}$ on $E(A_f)$.
\end{theorem}

\begin{proof}
Observe that $\pstab{\Aut(A_f^*)}{\bo A_f^*}$ is embedded in $\Aut(X)$, and $\pstab{\Aut(A_f)}{\bo A_f}$ is embedded in $\Aut(X_i)$.
Let $\{u,v\}={\bo A}={\bo A^*}$, where $u$ is not necessarily distinct from $v$.
Now the statement follows from Theorem~\ref{thm:expanded_inhomogeneous_product} by setting $X=A_f^*$ and $X_{i+1}=A_f$, with the two vertices $u$ and $v$ colored by different colors.
The coloring implies $\Aut(X_{i+1})=\pstab{\Aut(A_f)}{\bo A_f}$ and  $\Aut(X)=\pstab{\Aut(A_f^*)}{\bo A_f^*}$. With this identification in mind, the statement follows.
\end{proof}

\section{Spherical groups and equivariance of point-orbits}
\label{sec:spherical}

We call a finite subgroup of the orthogonal group $O(3)$ of $3\times 3$ orthogonal matrices a \emph{spherical group}.
Abstractly, these groups include infinite families $\gC_n$, $\gD_n$, $\gC_n\times\gC_2$, and $\gD_n\times\gC_2$, the symmetry groups of the platonic solids, namely, $\gA_5\times\gC_2$, $\gS_4\times\gC_2$, and $\gS_4$, and all the subgroups of these groups.
Recall that every  $3\times 3$ orthogonal matrix is uniquely determined by its action on the sphere $S^2 = \{x\in\mathbb{R}^3 : \|x\| = 1\}$.
Hence, there is a one-to-one correspondance between the elements of $O(3)$ and the \emph{isometries} of $S^2$.
There are $14$  possible types of spherical groups identified by the Conway's notation~\cite{conway2016symmetries}; see Table~\ref{fig:table}.

There are three basic types of isometries of $S^2$, namely rotations, reflections, and the antipodal involution, and every isometry is a composition of these.
A \emph{rotation} is an orientation-preserving isometry which fixes a pair of antipodal points, a \emph{reflection} is an orienation-reversing isometry fixing a great circle, and the \emph{antipodal involution} is an orientation-reversing isomotery swapping the pairs of antipodal points.

\begin{table}[t!]
\centering
\includegraphics{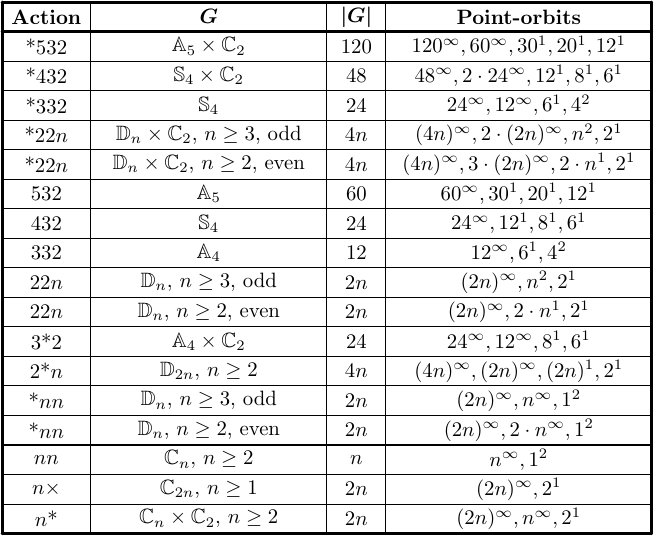}
\caption{
For each (parametrized) spherical group $G$, the first entry defines the action on the sphere, the second entry determines $G$ as an abstract group, the third entry gives the order of $G$.
The fourth entry gives the numbers of equivariance classes of
point-orbits of $G$.  By $c \cdot a^b$, we denote that there are $c$ equivariance
classes of point-orbits of size $a$ and each equivariance class consists of $b$ orbits.}
\label{fig:table}
\end{table}

\heading{Point-orbits.}
Let $g$ be an isometry of $S$ with a fixed point.
If $g$ is orientation-preserving, then $g$ is a \emph{rotation}, if $g$ is orienation-reversing, then $g$ is a \emph{reflection}.
It follows that a point-stabilizer in a spherical group is either cyclic or  a dihedral group.

For a spherical group $G$ and a point $x$ of the sphere, we distinguish four types of $O = [x]_G$ according to the structure of the point-stabilizers:
\begin{packed_itemize}
\item
\emph{regular}: $|O| = |G|$, the stabilizers of points in $O$ are trivial,
\item
\emph{regular-reflexive}: $|O| = |G|/2$ and the stabilizer of each point in $O$ is generated by a unique reflection,
\item
\emph{singular}: $|O| = |G|/m$, $m\geq 2$, and the stabilizer of each point in $O$ is generated by a rotation of order $m$,
\item
\emph{singular-reflexive}: $|O| = |G|/2m$, $m\geq 2$, and the stabilizer of each point in $O$ is generated by a rotation and reflection generating the dihedral group of order $2m$.
\end{packed_itemize}

With each spherical group $G$ there is an associated \emph{quotient orbifold} $\calO = S/G$ whose
points represent orbits of the action. In what follows, we denote by $\omega \colon S \to \calO$ the
natural projection $x \mapsto [x]_{G}$. Topologically, $\calO$ is a either a sphere, or a disk, or
the projective plane. The interior points of $\calO$ correspond to regular orbits and to singular
orbits. The boundary point of $\calO$ (if the quotient orbifold is topologically a disk) represents
a regular-reflexive, or a singular-reflexive orbit. There can be only a finite number of singular
and singular-reflexive orbits giving rise to finitely many \emph{branch points} on $\calO$. A branch-point $x$ is of index $m>1$, if $\omega^{-1}(x)$ is a singular orbit of size $|G|/m$, or if  $\omega^{-1}(x)$ is a reflexive singular orbit of size $|G|/2m$.  
The set
of all branch points will be denoted by $\calB$.

\heading{Regular orbits.}
For regular orbits the analysis is simple and follows trivially from Lemma~\ref{lem:trivial_orbits_equivariant}.

\begin{lemma}
\label{lem:regular_orbits}
All regular point-orbits in an action of $G$ on $S^2$ are equivariant.
\end{lemma}

\heading{Regular-reflexive orbits.}
If the set of regular-reflexive orbits is non-empty, then $\calO$ is topologically a disk and all the reflexive orbits are projected on the boundary.
From the classification of spherical groups, see Table~\ref{fig:table}, it follows that $|\calB \cap \bo \calO| \leq 3$.
Let $O_1$ and $O_2$ be two regular-reflexive orbits which are projected by $\omega$ to $\bar{x}$ and $\bar{y}$ in $\partial \calO$, respectively.
There are several cases.

\begin{packed_itemize}
\item
\emph{Case 1}: The points $\bar{x}$ and $\bar{y}$ are in the same connected component of
$\partial\calO\setminus \calB$.
\item
\emph{Case 2}: The points $\bar{x}$ and $\bar{y}$ are in different connected components of
$\bo\calO\setminus \calB$ incident to a branch point of an odd index.
\item
\emph{Case 3}: The points $\bar{x}$ and $\bar{y}$ are in different connected components of
$\bo\calO\setminus \calB$ and $\bar{x}$ is in a component separated by branch points of an even index.
\end{packed_itemize}

\begin{lemma}
\label{lem:reg_reflex_pos}
The regular-reflexive orbits $O_1 = \omega^{-1}(\bar{x})$ and $O_2 = \omega^{-1}(\bar{y})$ are equivariant in the Case 1 and 2.
\end{lemma}

\begin{proof}
In Case 1, there exists a path
$Q\subseteq \partial\calO$ joining $\bar{x}$ and $\bar{y}$ and avoiding the branch points (if there
are any). Then the $|G|/2$ lifts of $Q$ give a matching of $O_1$ and $O_2$, determining the equivariance
bijection.

Consider the induced action of $G$ on the connected components of $\omega^{-1}(\bo\calO\setminus\calB)$, which we call segments.
For $x \in S^2$, let $s_x$ be the segment containing the point $x$.
Then two orbits $[x]_G$ and $[y]_G$ are equivariant if and only if $[s_x]_G$ and $[s_y]_G$ are equivariant.
This observation reduces the identification of the equivariance classes to the problem of determining equivariance classes of orbits of the action of $G$ on the segments.
By Lemma~\ref{lem:equivariant_action_characterization}, $[s_x]_G$ and $[s_y]_G$ are equivariant if and only if $G_{s_x}$ and $G_{s_y}$ are conjugate in $G$.
We may assume that $s_x$ and $s_y$ are incident to $b \in \omega^{-1}(\bar b)$, for some $\bar b \in \calB$.
Since the action of $G$ is transitive on $\omega^{-1}(\bar b)$, the groups $G_{s_x}$ and $G_{s_y}$ are conjugate in $G$ if they are conjugate in $G_b \cong \gD_m$, where $m$ is the branch-index of $\bar b$.
Using Lemma~\ref{lem:dih_actions}, this concludes the proof for the Case 2.
\end{proof}

\begin{lemma}
\label{lem:reg_reflex_neg}
The regular-reflexive orbits $O_1 = \omega^{-1}(\bar{x})$ and $O_2 = \omega^{-1}(\bar{y})$ are non-equivariant in the Case 3.
\end{lemma}

\begin{proof}
Let $b_1$ and $b_2$ be the branch points of even index separating the segment containing $\bar{x}$. By the
classification of spherical groups, $G$ is one of the groups *432, *22n, and *nn ($n$ is
even).

First, we deal with the case *nn.  The boundary of $\calO$ lifts into an embedding $D_{2n}$ of a
$2n$-dipole on the sphere, where the two branch points of index $n$ lift to the two $2n$-valent
vertices of the dipole, respectively. Without a loss of generality, we can assume that the vertices
of the dipole correspond to the north and south pole of the sphere. The embedding determines a
cyclic ordering $e_1, e_2, \dots, e_{2n}$ of the edges of $D_{2n}$. We may assume that in this
ordering, every odd edge $e_i$ contains the unique preimage $x_i$ of $\bar{x}$ and every even edge
$e_j$ contains the unique
preimage $y_j$ of $\bar{y}$.

Let $\tau \in G$ be the reflection symmetry that fixes point-wise the great circle $C$ consisting of
$e_1$ and $e_{n+1}$. The circle $C$ separates the sphere $S$ into two open hemispheres $S^{+}$ and
$S^{-}$ transposed by $\tau$. Since $n$ is even, $e_1$ and $e_{n+1}$ contain two preimages $x_1$ and
$x_{n+1}$ of $\bar{x}$, respectively. For contradiction, assume that there exists an equivariance
bijection $f$ matching $O_1$ to $O_2$. In particular, this means that $f(x_1) = y$, where $y \notin
C$ is some preimage of $\bar{y}$. Without a loss of generality we can assume $y \in S^{+}$. The
mapping $f$ cannot be an equivariance bijection, since $f(\tau(x_1)) = f(x_1) = y \in S^{+}$ and
$\tau(f(x_1)) = \tau(y) \in S^{-}$.

Now we deal with the cases *432, *22n ($n$ is even). For each such action of $G$, the boundary of the
respective orbifold gives a barycentric subdivision of the cube, of an $n$-dipole. We
give the proof for *22n, the other case can be handled analogously. Similarly as above, $\calO$
lifts to a triangulation of the sphere determined by a $2n$-sided bi-pyramid with the two
$2n$-valent vertices located at the north and south pole. The embedding determines a natural cyclic
ordering $e_1,e_2,\dots, e_{2n}$ of the edges emanating from the vertex in the north pole and the
cyclic ordering $e_1',e_2',\dots, e_{2n}'$ at the south pole.  Similarly as in the previous
paragraph, every odd edge contains a unique preimage of $\bar{x}$ and every even edge contains a
preimage of $\bar{y}$. We may assume that in this ordering, the odd edges $e_i$ and $e_i'$ contain
the unique preimages $x_i$ and $x_i'$ of $\bar{x}$, respectively, and the even edges $e_j$ and
$e_j'$ contain the preimages $y_j$ and $y_j'$  of $\bar{y}$.

Let $\tau \in G$ be the reflection symmetry that fixes point-wise the great circle $C$ consisting of
$e_1$, $e_{n+1}$, $e_{n+1}'$, and $e_1'$. The circle $C$ separates the sphere $S$ into two open
hemispheres $S^{+}$ and $S^{-}$ transposed by $\tau$. Since $n$ is even, $e_1$, $e_{n+1}$,
$e_{n+1}'$, and $e_1'$ contain
preimages $x_1$, $x_{n+1}$, $x_{n+1}'$, and $x_1'$ of $\bar{x}$, respectively. In particular, $C
\cap O_2 = \emptyset$. For contradiction, assume that there
exists and equivariance bijection $f$ matching $O_1$ to $O_2$. This means that
$f(x_1) = y$, where $y \notin C$ is some preimage of $\bar{y}$. Without a loss of generality we can
assume $y \in S^{+}$. The mapping $f$ cannot be an equivariance bijection, since $f(\tau(x_1)) =
f(x_1) = y \in S^{+}$ and $\tau(f(x_1)) = \tau(y) \in S^{-}$.
\end{proof}

\heading{Singular and singular-reflexive orbits.}
The equivariance of singular and singular-reflexive orbits will be discussed for each spherical
group separately. It is sufficient to investigate the following
spherical groups: *332, *22n, 332, 22n, *nn, and nn.

\begin{lemma}
\label{lem:*332}
Let $G$ be a spherical group of type *332 or 332. Then the two orbits of size $4$ are equivariant.
\end{lemma}

\begin{proof}
The two orbits of size $4$ are represented by the vertices and centers of faces of the spherical
tetrahedron.
Observe that antipodal point to a point representing a vertex of the tetrahedron is the center of
the opposite face. This defines a matching between the vertices and the centers of faces of the
tetrahedron, giving the equivariance between the two singular orbits of size $4$.
\end{proof}

\begin{lemma}
\label{lem:*22n}
Let $G$ be a spherical group of type 22n or *22n. Then the following statements hold:
\begin{packed_itemize}
\item
If $n = 2$, then the three orbits of order $2$ are pairwise non-equivariant.
\item
If $n\geq 3$ and $n$ is odd, then the two orbits of order $n$ are equivariant.
\item
If $n\geq 4$ and $n$ is even, then the two orbits of order $n$ are non-equivariant.
\end{packed_itemize}
\end{lemma}

\begin{proof}
Let $G$ be of type $22n$.
Construct a geodesic triangle with vertices being the three branch points.
Then the triangle lifts to the $2n$-sided by-pyramid, where the orbits $O_1$ and $O_2$ form the vertices of the cycle $C$ of length $2n$ lying on the equator, thus the opposite points of $C$ are antipodal.
Now the statement follows from Lemma~\ref{lem:dih_actions} if $n > 2$.
If $n = 2$, then the by-piramid is the octahedron, where the three orbits of size two are formed by the three pairs of antipodal vertices.
The stabilizers of these orbits correspond to the three non-trivial subgroups of $G \cong \gC_2^2$.
Since $\gC_2^2$ is abelian, all conjugacy classes of subgroups are trivial, and the orbits are non-equivariant by Lemma~\ref{lem:equivariant_action_characterization}.

Let $G$ be of type $*22n$.
Then the orbifold $\calO$ lifts to the $2n$-sided by-pyramid, where the orbits $O_1$ and $O_2$ form the vertices of the cycle $C$ of length $2n$ lying on the equator, thus the opposite points of $C$ are antipodal.
Observe that $G$ acts on $C$ as the dihedral group $\gD_n$.
The statement follows from Lemma~\ref{lem:dih_actions} if $n > 2$.
If $n = 2$, then the by-piramid is the octahedron, where the three orbits of size two are formed by the three pairs of antipodal vertices.
The stabilizers of these orbits correspond to three distinct subgroups of $G \cong \gC_2^3$.
Since $\gC_2^3$ is abelian, all conjugacy classes of subgroups are trivial, and the orbits are non-equivariant by Lemma~\ref{lem:equivariant_action_characterization}.
\end{proof}

\begin{lemma}
\label{lem:*nn_nn}
Let $G$ be a spherical group of type *nn or of type nn. Then the two singular orbits of order $1$ are equivariant.
\end{lemma}

\begin{proof}
Follows trivially from Lemma~\ref{lem:trivial_orbits_equivariant}
\end{proof}

Using Lemmas~\ref{lem:regular_orbits}--\ref{lem:*nn_nn} we get the following theorem. 

\begin{theorem}
\label{thm:point}
The equivariance classes of point-orbits are described in the fourth column of Table~\ref{fig:table}.
\end{theorem}

\begin{figure}
\centering
\includegraphics{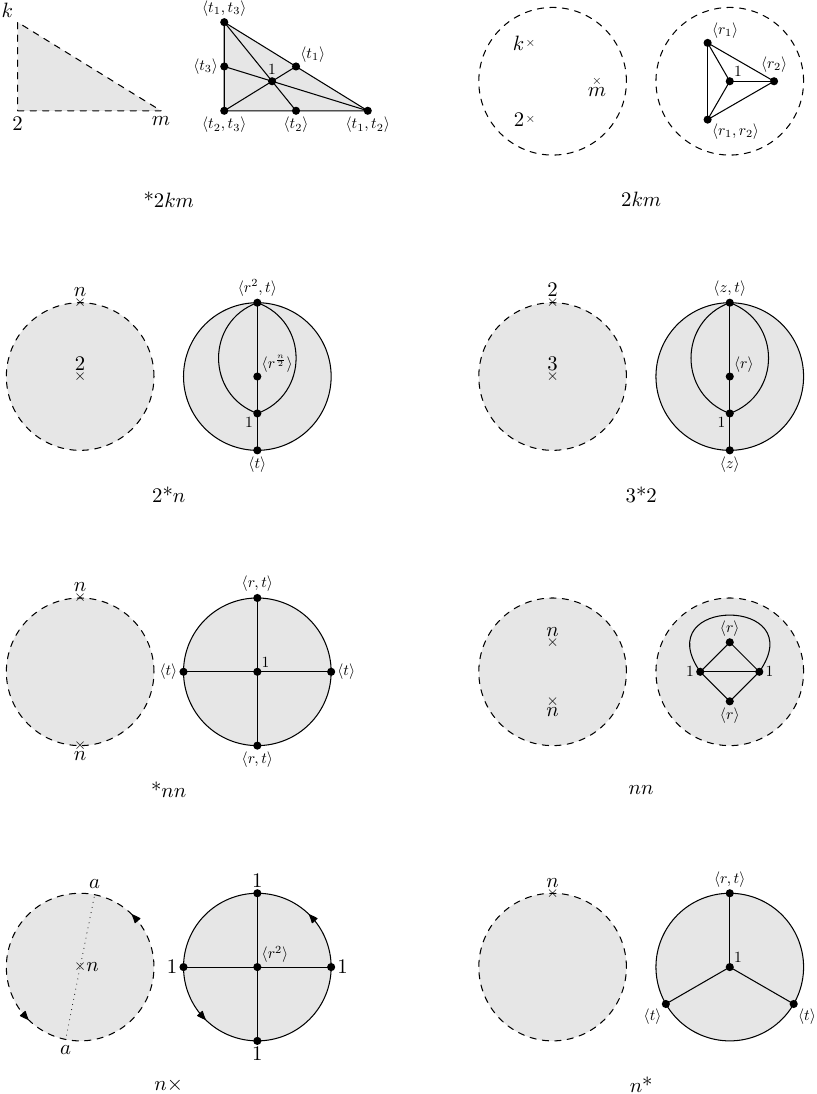}
\label{fig:orbifolds}
\caption{Quotient orbifolds and construction of polyhedra containing all types of orbits.}
\end{figure}

\heading{Constructions of polyhedra.}
In Figure~\ref{fig:orbifolds}, for each spherical group $G$ the quotient orbifold $O = S^2/G$ is depicted together with a graph $X$ embedded in $O$.
Each embedding lifts along the projection $S^2 \to O$ onto the $2$-skelenton of a polyhedron $P$.
The polyhedron $P$ satisfies the following properties:
\begin{packed_itemize}
\item
$G$ is a subgroup of $\Aut(P)$, and
\item
for every equivariance class of point-orbits of $G$ there exists a vertex-orbit in $P$ representing the clases.
\end{packed_itemize}

In each quotient graph $X$ in Figure~\ref{fig:orbifolds}, the vertices are labeled by subgroups of $G$.
For each vertex $v$ of $X$, the corresponding subgroup $H\leq G$ determines the action of $G$ on the fiber of $v$; see Lemma~\ref{lem:equivariant_action_characterization}.
In particular, we use the following presentations of the spherical groups:
\begin{packed_itemize}
\item 
$\mbox{*}2km$: $\langle t_1, t_2, t_3 \mid t_1^2 = t_2^2 = t_3^2 = (t_1t_3)^k = (t_1t_2)^m = (t_2t_3)^2 = 1\rangle$,
\item
$2km$: $\langle r_1, r_2 \mid r_1^k = r_2^m = (r_1r_2)^2 = 1\rangle$,
\item
$2\mbox{*}n$: $\langle r, t\mid r^{2n} = (rt)^2 = 1\rangle$,
\item
$3\mbox{*}2$: $\langle r, t \mid r^3 = t^2 = (rt)^3 = 1\rangle \times \langle z \mid z^2 =  1\rangle$,
\item
$\mbox{*}nn$: $\langle r,t \mid r^n = (rt)^2 = 1\rangle$,
\item
$nn$: $\langle r \mid r^n = 1\rangle$,
\item
$n\times$, $n$ even:  $\langle r \mid r^n = 1\rangle$,
\item
$n\mbox{*}$: $\langle r,t \mid r^n = (rt)^2 = 1\rangle$.
\end{packed_itemize}
Note that the above description of actions on the vertex-orbits does not require reference to spherical groups.

The polyhedrality of the lifted maps can be checked in each case  directly. For instance, it is not difficult to identify
the lifted maps in case $nn$, see Figure~\ref{fig:orbifolds}. These are the 2-skeletons
of $2n$-sided bipyramids. In case $*2km$ the lifted maps arise by a local operation applied to the barycentric subdivision of the $2$-skeletons of the five Platonic solids and to cycles (if one of the $k$ and $m$ is equal $2$). Such operations on maps were investigated in~\cite{heidi},
where a criterion implying polyhedrality is derived, in particular, the lifted maps in case $*2km$ are polyhedral.

\section{Automorphism groups of 3-connected planar graphs}
\label{sec:3connected}

In this section, we recall the classification of automorphism groups of $3$-connected planar graphs at the abstract level, and also at the level of actions on the sphere.
The classification is based on Whitney's theorem~\cite{whitney} stating that a 3-connected planar graph has a unique embedding into the sphere, and on Mani's theorem~\cite{mani1971automorphismen} establishing that the unique embedding can be realized on the sphere such that all automorphisms of the graph extend to isometries of the underlying sphere.

\heading{Maps on surfaces.}
A \emph{map} $\calM$ is a 2-cell decomposition of a compact connected surface $\Sigma$.
A map is usually defined by a 2-cell embedding of a connected graph $i \colon X \hookrightarrow \Sigma$.
The connectivity components of $\Sigma \setminus i(X)$ are called faces of $\calM$.  An automorphism of a map is an automorphism of the
graph preserving the incidences between the vertices, edges, and faces. If $\calM$ is a spherical map, then $\Aut(\calM)$ is one of
the spherical groups and with the exception of paths and cycles, it is a subgroup of $\Aut(X)$. As a
consequence of Whitney's theorem~\cite{whitney} we have the following.

\begin{theorem}\label{lem:aut_3_conn}
Let $\calM$ be the map given by the unique $2$-cell embedding of a $3$-connected graph $X$ into the
sphere. Then $\Aut(X) = \Aut(\calM)$.
\end{theorem}

By Mani's theorem~\cite{mani1971automorphismen}, there exists a polyhedron $P$, whose $1$-skeleton is isomorphic to $X$, such that the group of isometries of $P$ coincides with $\Aut(X)$.
Also, the polyhedron $P$ can be placed in the interior of
a sphere and projected onto it, so that each isometry of $P$ corresponds to an isometry of the sphere.
Therefore, every automorphism in $\Aut(X)$ can be viewed as an isometry of the sphere.
In particular, we have:

\begin{theorem} \label{thm:mani}
Let $X$ be a $3$-connected planar graph. Then $\Aut(X)$ is isomorphic to one of the spherical groups.
\end{theorem}

We recall some basic definitions from geometry~\cite{stillwell,vca}. An automorphism of a
3-connected planar graph $G$ is called \emph{orientation preserving}, if the respective isometry
preserves the global orientation of the sphere. It is called \emph{orientation reversing} if it
changes the global orientation of the sphere. A subgroup of $\Aut(X)$ is called \emph{orientation
preserving} if all its automorphisms are orientation preserving, and \emph{orientation reversing}
otherwise. Since the composition of orientations reversing automorphisms is orientations preserving, every orientation reversing subgroup contains an orientation preserving
subgroup of index 2.

\heading{Vertex- and edge-stabilizers.}
Let $X$ be a $3$-connected planar graph and let $u \in V(X)$.
The stabilizer of $u$ in $\Aut(X)$ is a subgroup of a dihedral group and it has
the following description in the language of isometries. If $\Aut(X)_u \cong \gC_n$, for $n \ge 3$,
it is generated by a rotation of order $n$ that fixes $u$ and the opposite point of the
sphere, and fixing no other point of the sphere. The opposite point of the sphere may be another
vertex or a center of a face.  If $\Aut(X)_u \cong \gD_n$, for $n \ge 2$, it consists of rotations
fixing $u$ and the opposite point of the sphere and reflections fixing a great circle passing
through $u$ and the opposite point. Each reflection always fixes either a center of some edge, or
another vertex. When $\Aut(X)_u \cong \gD_1 \cong \gC_2$, it is generated either by a $180^\circ$
rotation, or by a reflection.

Let $e \in E(X)$. The stabilizer of $e$ in $\Aut(X)$ is a subgroup of $\gC_2^2$. When $\Aut(X)_e
\cong \gC_2^2$, it contains the following three non-trivial isometries. First, the $180^\circ$
rotation around the center of $e$ and the opposite point of the sphere that is a vertex, center of
an edge, or center of an even face. Next, two reflections perpendicular to each other which fix the center of  
 $e$, the opposite point of the sphere. When $\Aut(X)_e \cong \gC_2$, it is generated by
only one of these three isometries.

We summarize the above discussion in the following lemma.

\begin{lemma}
\label{lem:vertex_edge_stabilizers}
Let $X$ be a $3$-connected planar graph, let $u$ be a vertex of degree $n$, and let $e$ be an edge.
Then $\Aut(X)_u$ is isomorphic to a subgroup of $\gD_n$, and $\Aut(X)_e$ is isomorphic to a subgroup of $\gC_2^2$.
\end{lemma}

\heading{Automorphism groups of planar primitive graphs.}

\begin{lemma} \label{lem:planar_primitive_graph}
The automorphism group $\Aut(X)$ of a planar primitive graph $X$ is a spherical group.
\end{lemma}

\begin{proof}
By Lemma~\ref{lem:structure_primitive}, the essence $\overline{X}$ of $X$ is 3-connected, or an $n$-cycle, or $K_2$, or $K_1$.
If $\overline{X}$ is 3-connected, then $\Aut(\overline{X})$ is a spherical group by Theorem~\ref{thm:mani}.
If $\overline{X}$ is $K_1$, $K_2$ or $C_n$ with attached single pendant edges, then $\Aut(\overline{X})$ is a subgroup of $\gC_2$, or of $\gD_n$.
Clearly, $\Aut(X)$ is a subgroup of $\Aut(\overline{X})$.
Since the family of spherical groups is closed under taking subgroups, we are done.
\end{proof}

\heading{Stabilizers of planar atoms.}

\begin{lemma}\label{lem:planar_atom_aut_groups}
Let $A$ be a planar atom.
\begin{packed_enum}
\item[(a)] If $A$ is a star atom, then $\sstab{\Aut(A)}{\bo A} = \pstab{\Aut(A)}{\bo A}$ which is a direct product of
symmetric groups. 
\item[(b)] If $A$ is a block atom, then $\sstab{\Aut(A)}{\bo A} = \pstab{\Aut(A)}{\bo A}$ and it is a subgroup of
a dihedral group.
\item[(c)] If $A$ is a proper atom, then $\sstab{\Aut(A)}{\bo A}$ is a subgroup of $\gC_2^2$ and $\pstab{\Aut(A)}{\bo A}$ is
a subgroup of $\gC_2$.
\item[(d)] If $A$ is a dipole, then $\pstab{\Aut(A)}{\bo A}$ is a direct product of symmetric groups. If $A$ is
symmetric, then $\sstab{\Aut(A)}{\bo A} = \pstab{\Aut(A)}{\bo A} \times \gC_2$. If $A$ is asymmetric, then $\sstab{\Aut(A)}{\bo A} =
\pstab{\Aut(A)}{\bo A}$.
\end{packed_enum}
\end{lemma}

\begin{proof}
(a) Since $|\bo A|=1$, we have $\sstab{\Aut(A)}{\bo A} = \pstab{\Aut(A)}{\bo A}$.  
Since the edges of each color class of the star atom $A$ can be arbitrarily and independently permuted, $\pstab{\Aut(A)}{\bo A}$ is a direct product of symmetric groups.

(b) Similarly as in Case (a), $|\bo A|=1$, and we have $\sstab{\Aut(A)}{\bo A} = \pstab{\Aut(A)}{\bo A}$.
Let $B = \overline{A}$.
It follows that $\sstab{\Aut(A)}{\bo A}\leq\sstab{\Aut(B)}{\bo B}$. By
Lemma~\ref{lem:structure_block}, either $B$ is a cycle, $K_2$, or a 3-connected
planar graph. In the first two cases, $\sstab{\Aut(B)}{\bo B}$ is a subgroup of $\gC_2$, while in the last case, it is the
stabilizer of a vertex in a 3-connected planar graph which is a subgroup of $\gD_n$, where $n$ is the degree
of the articulation separating $A$.

(c) Let $A$ be a proper atom with $\bo A = \{u,v\}$.
Let $B = \overline{A}$, so $\sstab{\Aut(A)}{\bo A} \leq \sstab{\Aut(B)}{\bo B}$ and $\pstab{\Aut(A)}{\bo A} \leq \pstab{\Aut(B)}{\bo B}$.
Clearly, $\sstab{\Aut(B)}{\bo B} = \sstab{\Aut(B^+)}{\bo B^+}$, where $B^+=B+ uv$.
By Lemma~\ref{lem:structure_proper}, $B^+$ is either a cycle, or a 3-connected  graph. In the
former case, $\Aut(B^+)$ is a subgroup of $\gC_2$ and $\pstab{\Aut(B^+)}{\bo B^+}$ is trivial.
In the latter case, we claim that $B^+$ is planar.
First, by the definition of a proper atom, $B$ is a subgraph of a block in the planar graph $X$.
It follows that $X\setminus \int B$ is a connected plane graph, in particular, the boundary vertices $u$ and $v$ are connected by a path in $X\setminus \int B$.
Thus, in the embedding of $B$ into the sphere the vertices $u$ and $v$ appear in the boundary of the same face, i.e.,  $B^+$ is planar. 
By Lemma~\ref{lem:vertex_edge_stabilizers}, $\Aut(B^+)_e$ of $e = uv$ is isomorphic to a subgroup of $\gC_2^2$.
It follows that  $\sstab{\Aut(A)}{\bo A}$ is a subgroup of $\gC_2^2$ and $\pstab{\Aut(A)}{\bo A}$ is a subgroup of $\gC_2$.

(d) Let  $\{u,v\}$ be the vertex set of the dipole $A$.
If $A$ is asymmetric, we have $\sstab{\Aut(A)}{\bo A} = \pstab{\Aut(A)}{\bo A}$ which is a direct product of symmetric groups. 
If $A$ is symmetric, there is an involution $t \in \sstab{\Aut(A)}{\bo A}$ which swaps $u$ with $v$ and fixes all the edges.
Clearly, we have $\langle t\rangle \cap \pstab{\Aut(A)}{\bo A} = \{id\}$ and $t$ centralizes $\pstab{\Aut(A)}{\bo A}$.
\end{proof}

\begin{figure}[t!]
\centering
\includegraphics{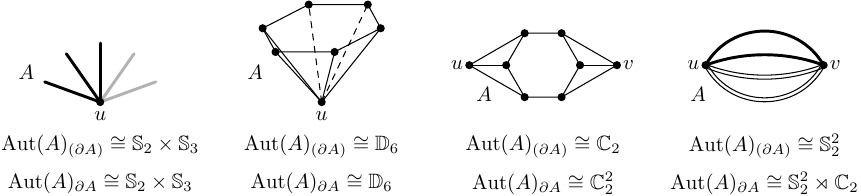}
\caption{An atom $A$ together with its groups $\pstab{\Aut(A)}{\bo A}$ and $\sstab{\Aut(A)}{\bo A}$.  From left to right, a star
block atom, a non-star block atom, a proper atom, and a dipole.}
\label{fig:automorphisms_of_atoms}
\end{figure}

By the parts (c) and (d) of Lemma~\ref{lem:planar_atom_aut_groups}, we obtain the following statement.

\begin{lemma} \label{lem:symmetric_proper_atom_involutions}
Every planar symmetric proper atom and dipole atom is centrally symmetric.
\end{lemma}

\section{Jordan-like characterization}
\label{sec:characterization}

In this section, we prove the main results of this paper: a complete recursive description of the automorphism groups of connected planar graphs.
Note that in combination with Theorem~\ref{thm:aut_disconnected}, this gives the description of the automorphism groups of all planar graphs.
The first step is to determine the abstract groups that can be realized as vertex-stabilizers in  Section~\ref{sec:fix_groups}.
The second step is to determine how these groups are composed with spherical groups in Section~\ref{sec:composition}.

\subsection{Vertex-stabilizers of planar graphs} \label{sec:fix_groups}

The aim of this section is to analyze the following class of abstract groups
$$\Stab(\planar) = \{G \cong \Aut(X)_u : X \in \planar, u \in V(X)\}$$
The following lemma relates $\Stab(\planar)$ together with pointwise stabilizers of boundaries of expanded atoms, allowing us to apply the theory built in Section~\ref{sec:reduction}.

\begin{lemma}
\label{lem:stab_fix}
The class $\Stab(\planar) = \Fix(\planar)$, where
$$\Fix(\planar) = \{G\cong \pstab{\Aut(X)}{\bo X} : X \in\planar, |E(X^R)| = 1, \bo X = V(X^R)\}.$$
\end{lemma}

\begin{proof}
Let, $G\in \Fix(\planar)$ and let $X$ be a graph certifying that.
Let $v \in \bo X$.
If $\bo X = \{v\}$, then $G \cong \pstab{\Aut(X)}{\bo X} = \Aut(X)_v \in \Stab(\planar)$.
If $\bo X = \{u, v\}$, then we form the graph $Y$ by indentifying the vertex $u$ with an end-vertex of a path $P$, where $|E(P)| > |E(X)|$.
Then $G \cong \pstab{\Aut(X)}{\bo X} \cong \Aut(Y)_v \in \Stab(\planar)$.

Let $G \in \Stab(\planar)$ and let $X$ and $v \in V(X)$ be a graph and a vertex certifying that.
Let $Y$ be the graph formed by identifying $v$ with an end-vertex of a path $P$, where $|E(P)| > |E(X)|$.
Now, the graph $Y^R$ is either a two-vertex graph, or $K_1$ with a single pendant edge attached, depending on the parity of $|E(P)|$.
The latter case can be alsways ensured by suitably adjusting the parity of $|E(P)|$.
We have $|E(Y^R)| = 1$ and $G \cong \Aut(X)_u \cong \pstab{\Aut(Y)}{\bo Y} \in \Fix(\planar)$.
\end{proof}

The next theorem gives a recursive characterization of $\Stab(\planar)$.

\begin{theorem} \label{thm:stabilizers}
The class $\Stab(\planar) = \calF$, where $\calF$ is defined inductively as follows:
\begin{packed_itemize}
\item[(a)] $\{1\} \in \calF$.
\item[(b)] If $G_1,G_2 \in \calF$, then $G_1 \times G_2 \in \calF$.
\item[(c)] If $G \in \calF$, then $G \wr \gS_n \in \calF$.
\item[(d)] If $G \in \calF$, then $G \wr \gC_n \in \calF$.
\item[(e)] If $G_1,G_2\in \calF$, then $(G_1, G_2) \gwr \gD_n \in \calF$, for $n\geq 3$ odd, where $|\Omega_1| = 2n$ and $|\Omega_2| = n$.
\item[(f)] If $G_1,G_2,G_3\in \calF$, then $(G_1, G_2, G_3) \gwr \gD_n \in \calF$, for $n\geq 2$ even,
where $|\Omega_1| = 2n$ and $|\Omega_2| = |\Omega_3| = n$, and $\gD_n$ acts on $\Omega_1$ regularly, and acts on $\Omega_2$ and $\Omega_3$ as on the vertices and edges of the regular $n$-gon, respectively.
\end{packed_itemize}
\end{theorem}

We split the proof into the following two lemmas.

\begin{lemma}
Each group in $\Fix(\planar)$ is isomorphic to a group in $\calF$.
\end{lemma}
\begin{proof}
Let $X = X_0,X_1,\dots,X_r$ be the reduction series of a planar graph $X$.
We proceed by an induction on $i = 1,\dots,r$.
If $i = 1$, then $A_e^* = A_e$, for any $e \in E(X_1)$, is an atom of $X = X_0$.
By Lemma~\ref{lem:planar_atom_aut_groups}, we have $\pstab{\Aut(A_e)}{\bo A_e} \in \calF$.

If $i > 1$, we fix an edge $f \in E(X_i)$ and set $A = A_f$.
By Theorem~\ref{thm:semidirect_product},
$$\pstab{\Aut(A^*)}{\bo A^*} \cong \left (\prod_{e\in E(A)} \pstab{\Aut(A_e^*)}{\bo A_e^*} \right )\rtimes \pstab{\Aut(A)}{\bo A},$$
where the right hand side is the inhomogeneous wreath product defined by the action of $\pstab{\Aut(A)}{\bo A}$ on $E(A)$.
Recall that if an edge $e \in E(A)$ does not expand into an atom, then $\pstab{\Aut(A_e^*)}{\bo A_e^*}$ is trivial.

Let $\{\Omega_1,\dots,\Omega_m\}$ be the partition of $E(A)$ into the orbits of the action of $\pstab{\Aut(A)}{\bo A}$.
By induction hypothesis, for every $j \in \{1,\dots,m\}$ and for every $e \in \Omega_j$ and there is  $K_j \in \calF$ such that $K_j \cong \pstab{\Aut(A_{e}^*)}{\bo A_{e}^*}$.
Then
\begin{equation}\label{eq:pstab_orb}
\pstab{\Aut(A^*)}{\bo A^*} \cong G = \left (K_1^{\ell_1}\times \cdots \times K_m^{\ell_m} \right ) \rtimes S,
\end{equation}
where $\ell_i = |\Omega_i|$, and $S\cong \pstab{\Aut(A)}{\bo A}$ acts on $E(A)$.
We split the proof into several cases, depending on the type of $A$.

\emph{Case 1: $A$ is a star atom or a dipole atom.}
By Lemma~\ref{lem:planar_atom_aut_groups}, $S = \gS_{\ell_1}\times\cdots\times \gS_{\ell_m}$, where $\gS_{\ell_j}$ is isomorphic to the subgroup of $\pstab{\Aut(A)}{\bo A}$ fixing every edge not in the orbit $\Omega_j$.
For every $j$, $G$ contains a subgroup $G_j \cong K_j \wr \gS_{\ell_j} \in \calF$.
By definition, the subgroups $G_j$ and $G_k$, for $j \neq k$, are disjoint.
Moreover, if $g \in G_j$ and $g' \in G_k$ then $gg' = g'g$.
Therefore,
$$\pstab{\Aut(A^*)}{\bo A^*}\cong G \cong  K_1\wr\gS_{\ell_1} \times \cdots \times K_m\wr\gS_{\ell_m}.$$
The group on the right side belongs to $\calF$.

\emph{Case 2: $A$ is a proper atom.}
By Lemma~\ref{lem:planar_atom_aut_groups}, either $\pstab{\Aut(A)}{\bo A}$ is trivial, or $\pstab{\Aut(A)}{\bo A} \cong \gC_2$.
In the first case, there is nothing to prove since each $\pstab{\Aut(A_e^*)}{\bo A_e^*}$, for $e \in E(A)$, is isomorphic to a group from $\calF$ and so is their direct product.

In the second case, the induced action of $\gC_2$ on $E(A)$ has orbits $\{\Omega_1,\dots,\Omega_m\}$ of size at most $2$.
Let $\Omega_1,\dots,\Omega_t$, for some $0\leq t \leq m$, be the orbits of size $1$.
By substituting to~(\ref{eq:pstab_orb}), we get
$$G = \left(K_1\times\cdots\times K_t \times K_{t+1}^2\times\cdots\times K_{m}^2 \right) \rtimes \gC_2,$$
where $K_i \in \calF$, for $i=1,\dots,m$, by induction.
Denote $G_1 =K_1\times\cdots\times K_t$ and  $G_2 =  K_{t+1}\times\cdots\times K_{m}$.
By Lemma~\ref{lem:inhomogeneous_fixed}
$$G \cong G_1\times (K_{t+1},\dots,K_m)\gwr \gC_2.$$
By Lemma~\ref{lem:inhomogeneous_equivariant}
$$G \cong G_1\times G_2\wr \gC_2,$$
which belongs to $\calF$.

\emph{Case 3: $A$ is a block atom.}
Since $A$ is essentially $3$-connected (Lemma~\ref{lem:structure_block}), $\pstab{\Aut(A)}{\bo A}$ is a spherical group. 
By Lemma~\ref{lem:planar_atom_aut_groups}, $\pstab{\Aut(A)}{\bo A} \cong S$ is a subgroup of a dihedral group, in particular, either $S$ cyclic of type $nn$ or $S$ is dihedral of type $\mbox{*}nn$.
We distinguish three cases: $S \cong \gC_n$ with $n\geq 2$, $S \cong \gD_n$ with $n\geq 3$ and $n$ odd, $S \cong \gD_n$ with $n\geq 2$ and $n$ even.

Let $\pstab{\Aut(A)}{\bo A} \cong \gC_n$.
The induced action of $\gC_n$ on $E(A)$ has orbits $\{\Omega_1,\dots,\Omega_m\}$ of size $1$ or $n$; see Table~\ref{fig:table}.
Let $\Omega_1,\dots,\Omega_t$ be the orbits of size $1$, for some $0\leq t \leq 2$.
By substituting to~(\ref{eq:pstab_orb}), we get
$$G = \left(K_1\times\cdots\times K_t \times K_{t+1}^n\times\cdots\times K_{m}^n \right) \rtimes \gC_n,$$
where $K_i \in \calF$, for $i=1,\dots,m$, by induction.
Denote $G_1 =K_1\times\cdots\times K_t$ and  $G_2 =  K_{t+1}\times\cdots\times K_{m}$.
By Lemma~\ref{lem:inhomogeneous_fixed}
$$G \cong G_1\times (K_{t+1},\dots,K_m)\gwr \gC_n.$$
Since the action of $\gC_n$ is regular on each $\Omega_j$, for $j>t$, the actions of $\gC_n$ on $\Omega_j$ are all equivariant.
By Lemma~\ref{lem:inhomogeneous_equivariant}
$$G \cong G_1\times G_2\wr \gC_n,$$
which belongs to $\calF$.

Let $\pstab{\Aut(A)}{\bo A} \cong \gD_n$ with $n \geq 3$ and $n$ odd.
The induced action of $\gD_n$ on $E(A)$ has orbits $\{\Omega_1,\dots,\Omega_m\}$ of possible sizes $1$, $n$, and $2n$; see Table~\ref{fig:table}.
There are integers $0\leq t_1 \leq t_2 \leq m$, $t_1\leq 2$, such that the orbits $\Omega_1,\dots,\Omega_{t_1}$ are of size $1$, the orbits $\Omega_{t_1+1},\dots, \Omega_{t_2}$ are of size $n$, and the orbits $\Omega_{t_2+1},\dots,\Omega_m$ are of size $2n$.
By substituting to~(\ref{eq:pstab_orb}), we get
$$G = \left(K_1\times\cdots\times K_{t_1} \times K_{t_1+1}^n\times\cdots\times K_{t_2}^n \times K_{t_2+1}^{2n}\times\cdots\times K_{m}^{2n} \right) \rtimes \gD_n,$$
where $K_i \in \calF$, for $i=1,\dots,m$, by induction.
Denote $G_1 =K_1\times\cdots\times K_{t_1}$, $G_2 =  K_{t_1+1}\times\cdots\times K_{t_2}$, and $G_3 = K_{t_2+1}\times\cdots\times K_{m}$.
By Lemma~\ref{lem:inhomogeneous_fixed}
$$G \cong G_1\times (K_{t_1+1},\dots,K_m)\gwr \gD_n.$$
By Lemma~\ref{lem:dih_actions}, the actions of $\gD_n$ on all $\Omega_j$, for $t_1 < j \leq t_2$, are equivariant.
The actions of $\gD_n$ on all $\Omega_j$, for $j>t_2$, are all regular and therefore equivariant.
By Lemma~\ref{lem:inhomogeneous_equivariant}
$$G \cong G_1\times (G_2,G_3)\gwr \gD_n,$$
which belongs to $\calF$.

Let $\pstab{\Aut(A)}{\bo A} \cong \gD_n$ with $n \geq 2$ and $n$ even.
The induced action of $\gD_n$ on $E(A)$ has orbits $\{\Omega_1,\dots,\Omega_m\}$ of possible sizes $1$, $n$, and $2n$; see Table~\ref{fig:table}.
By Lemma~\ref{lem:dih_actions}, there are two possible non-equivariant actions of $\gD_n$ on $n$ points.
Thus, there are integers $0\leq t_1 \leq t_2 \leq t_3 \leq m$ such that the orbits $\Omega_1,\dots,\Omega_{t_1}$ are of size $1$, the orbits $\Omega_{t_1+1},\dots, \Omega_{t_2}$ are of size $n$, the orbits $\Omega_{t_2+1},\dots,\Omega_{t_3}$ are of size $n$, and the orbits $\Omega_{t_3+1},\dots,\Omega_m$ are of size $2n$, and the $\gD_n$-sets $\Omega_i$ and $\Omega_j$, for $t_1 < i \leq t_2$ and $t_2 < j \leq t_3$, are non-equivariant.
By substituting to~(\ref{eq:pstab_orb}), we get
$$G = \left(\prod_{i=1}^{t_1}K_i \times \prod_{i=t_1+1}^{t_2}K_i^n \times\prod_{i=t_2+1}^{t_3}K_i^n \times\prod_{i=t_3+1}^{m}K_i^{2n} \right) \rtimes \gD_n,$$
where $K_i \in \calF$, for $i=1,\dots,m$, by induction.
Denote $G_1 =K_1\times\cdots\times K_{t_1}$, $G_2 =  K_{t_1+1}\times\cdots\times K_{t_2}$, $G_3 = K_{t_2+1}\times\cdots\times K_{t_3}$, and $G_4 = K_{t_3+1}\times\cdots\times K_{m}$.
By Lemma~\ref{lem:inhomogeneous_fixed}
$$G \cong G_1\times (K_{t_1+1},\dots,K_m)\gwr \gD_n.$$
By Lemma~\ref{lem:inhomogeneous_equivariant}
$$G \cong G_1\times (G_2,G_3,G_4)\gwr \gD_n,$$
which belongs to $\calF$.
\end{proof}

\begin{figure}[t!]
\centering
\includegraphics{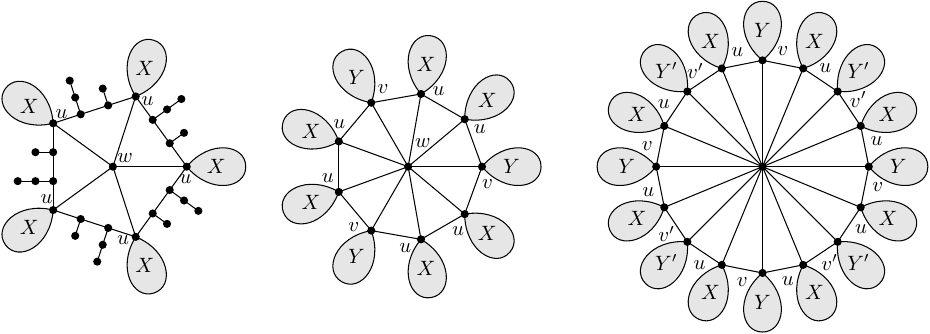}
\caption{(a) $\Aut(Z) \cong \Aut(X)\wr \gC_5$, (b) $\Aut(Z) \cong (\Aut(X), \Aut(Y))\gwr \gD_3$, (c) $\Aut(Z) \cong (\Aut(X), \Aut(Y), \Aut(Y'))\gwr \gD_4$}
\label{fig:wheels}
\end{figure}

\begin{lemma} \label{lem:fix_groups_construction}
Each group in $\calF$ is isomorphic to a group in $\Stab(\planar)$.
\end{lemma}

\begin{proof}
We prove the statement by induction on the number of operations needed to construct a group in $\calF$.
Clearly, the trivial group is in $\Stab(\planar)$.
Let $H \in \calF$ be nontrivial, then it is constructed by one of the operations (b)--(f), for which we prove the statement separately.

Let $H = G_1\times G_2 \in \calF$, for $G_1,G_2 \in \Stab(\planar)$.
By the induction hypothesis, there are $X$ and $Y$ such that $G_1\cong \Aut(X)_u$ and $G_2\cong \Aut(Y)_v$.
If $X \cong Y$, we modify $Y$ by attaching an end-vertex of a copy of $P_2$ to all vertices in $V(Y)\setminus \{v\}$.
We form a $Z$ with $w \in V(Z)$, where $w$ arises by identifying $u$ and $v$.
Now, $H = G_1\times G_2 \cong \Aut(Z)_w \in \Stab(\planar)$.

Let $H = G\wr \gS_n \in \calF$, for $G \in \Stab(\planar)$.
By the induction hypothesis, there is a graph $X$ and $u$ such that $G \cong \Aut(X)_u$.
We form a $Z$ with $w \in V(Z)$, where $w$ is arises by identifying $n$ copies of $X$ at the vertex $u$.
Now, $H = G \wr \gS_n\cong \Aut(Z)_w \in \Stab(\planar)$.

Let $H = G\wr \gC_n \in \calF$, for $G \in \Stab(\planar)$.
By the induction hypothesis, there is a graph $X$ and $u$ such that $G \cong \Aut(X)_u$.
We form a graph $Z$ as follows.
We take an $n$-wheel $W_n$ with the center $w$ and outer vertices labeled $\{0,\dots,n-1\}$.
We indentify the vertex $i$ with the vertex $u$ in a copy of $X$.
Finally, we replace every edge joining the vertices $i$ and $i+1 \mod n$ by the graph below.
\begin{center}
\includegraphics{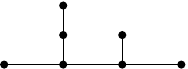}
\end{center}
Now, $H = G \wr \gC_n \cong \Aut(Z)_w \in \Stab(\planar)$.

Let $H = (G_1,G_2)\gwr \gD_n \in \calF$, for $G_1,G_2 \in \Stab(\planar)$ and $n\geq 3$ odd.
By the induction hypothesis, there are graphs $X$ and $Y$ such that $G_1 \cong \Aut(X)_u$ and $G_2 \cong \Aut(Y)_v$.
We form a graph $Z$ as follows.
We take an $3n$-wheel $W_{3n}$ with the center $w$ and outer vertices labeled $\{0,\dots,3n-1\}$.
For every $i = 0 \pmod 3$, we indentify the vertex $i$ with the vertex $v$ in a copy of $Y$.
For every $i \neq 0 \pmod 3$, we identify the vertex $i$ with the vertex $u$ in a copy of $X$.
Now, $H = (G_1,G_2) \gwr \gD_n\cong \Aut(Z)_w \in \Stab(\planar)$.

Let $H = (G_1,G_2, G_3)\gwr \gD_n \in \calF$, for $G_1,G_2, G_3 \in \Stab(\planar)$ and  $n\geq 2$ even.
By the induction hypothesis, there are graphs $X$, $Y$, and $Y'$ such that $G_1 \cong \Aut(X)_u$, $G_2 \cong \Aut(Y)_v$, and $G_3\cong \Aut(Y')_{v'}$.
We form a graph $Z$ as follows.
We take an $4n$-wheel $W_{4n}$ with the center $w$ and outer vertices labeled $\{0,\dots,4n-1\}$.
For every $i = 1 \pmod 2$, we indentify the vertex $i$ with the vertex $u$ in a copy of $X$.
For every $i = 0 \pmod 4$, we identify the vertex $i$ with the vertex $v$ in a copy of $Y$.
For every remaining $i$, we identify the vertex $i$ with the vertex $v'$ in a copy of $Y'$.
Now, $H = (G_1,G_2, G_3) \gwr \gD_n\cong \Aut(Z)_w \in \Stab(\planar)$.
\end{proof}

\subsection{Composition of spherical groups with vertex-stabilizers}\label{sec:composition}

Finally, can state our main theorem.

\begin{theorem}
\label{thm:composition}
We have $\Aut(\planar) = \calP$, where $\calP$ is defined as follows:
if $G_1,\dots, G_m \in \Stab(\planar)$ and $S$ is a spherical group acting on $\Omega$ with pairwise non-equivariant orbits $\Omega_1,\dots,\Omega_m$, then $(G_1,\dots,G_m)\gwr S \in \calP$.
\end{theorem}

Recall that for every spherical group $S$, the possible equivariance classes of orbits of the action of $S$ on a set $\Omega$ is are described in Table~\ref{fig:table}.

\begin{proof}[Proof of Theorem~\ref{thm:composition}]
By Lemma~\ref{lem:stab_fix}, we have that $\Stab(\planar) = \Fix(\planar)$, thus we shall use them interchangebly.

$\Aut(\planar)\subseteq \calP$.
Let $X$ be a planar graph and let $S \cong \Aut(X^R)$.
By Theorem~\ref{thm:expanded_inhomogeneous_product}, $\Aut(X) \cong (K_1,\dots,K_\ell)\gwr_{\Omega'} S$, where $\Omega' = E(X^R)$, and $\ell$ is the number of edge-orbits in the action of $S$ on $\Omega'$, and $G_i \in \Fix(\planar)$, for $i=1,\dots,\ell$.
By Theorem~\ref{thm:mani}, $S$ is a spherical group.
By Lemma~\ref{lem:inhomogeneous_equivariant}, there are groups $G_1,\dots,G_m \in \Fix(\planar)$ such that $\Aut(X) \cong  (G_1,\dots,G_{m})\gwr_{\Omega} S$, where the $S$-sets $\Omega_i$, for $i = 1,\dots,m$, are pairwise non-equivariant.

$\Aut(\planar)\supseteq \calP$.
If $G_1,\dots, G_m \in \Stab(\planar)$ and $S$ is a spherical group acting on $\Omega$ with pairwise non-equivariant orbits $\Omega_1,\dots,\Omega_m$, then let $G = (G_1,\dots,G_m)\gwr S \in \calP$.
First, we construct a $3$-connected primitive planar graph $Y$, with $\Aut(Y)\cong S$, and with orbits $V_1,\dots,V_m \subseteq V(Y)$ such that the $(S,V_i)$ and $(S, \Omega_i)$ are equivariant.
Clearly, such a graph $Y$ is uniquely determined by its quotient $Y/S$ embedded in the corresponding orbifold.
For each spherical group $S$, such a quotient is depicted in Figure~\ref{fig:orbifolds}.
The reader may easily verify that the lift $Y$ is in each case $3$-connected.
It may happen that $S$ is isomorphic to a proper subgroup of $\Aut(Y)$, in this case we use a suitable coloring of edges to ensure $\Aut(Y)\cong S$.

Now, to construct $X$, we proceed as follows.
For $i=1,\dots,m$, we a choose graph $Y_i$ with $G_i \cong \pstab{\Aut(Y_i)}{\bo Y_i}$, and for every $v\in V_i$, we take a copy of $Y_i$ and identify the unique vertex of $\bo Y_i$ with $v$.
\end{proof}

\bibliographystyle{plain}
\bibliography{aut_planar}

\appendix

\section{Proofs of structural lemmas for atoms}
\label{app:lemmas}

To keep the paper self-contained, in this appendix, we include the proofs of several lemmas concerning the structural properties of atoms.

\begin{proof}[Proof of Lemma~\ref{lem:essence}]
If $\overline{X}$ has an articulation, then it contains a block part and also an atom, which is a contradiction.
Hence, $\overline{X}$ is $2$-connected.
Let $u$ be the vertex of minimum degree in $\overline{X}$.
If $\deg(u) = 1$, then $\overline{X} \cong K_2$, which is a contradiction.
If $\deg(u) = 2$, then either $\overline{X} \cong C_n$, or $u$ is an inner vertex of an induced path connecting two vertices $x$ and $y$ of degree at least $3$, or if $X$ is an extended proper atom, then $\{x,y\} = \bo X$.
In the former case the path is a proper atom, which is a contradiction, and in the latter case, $X$ is a cycle.
Finally, if $\deg(u) \geq 3$, then every $2$-cut is non-trivial.
It follows that $\overline{X}$ contains a proper part and therefore also an atom, which is a contradiction.
Thus, $\overline{X}$ has no $2$-cut and it is $3$-connected.
\end{proof}

\begin{proof}[Proof of Lemma~\ref{lem:structure_primitive}]
If $X$ has at least $3$ vertices, then by Lemma~\ref{lem:essence}, $\overline{X}$ is $3$-connected or a cycle.
Since $X$ cannot contain a pendant star, $X$ consists of $\overline{X}$ with possibly single pendant edges attached.
Such a graph is primitive since it has no parts.
If $X$ has at most $2$ vertices, there are exactly $5$ possible primitive graphs.
\end{proof}

\begin{proof}[Proof of Lemma~\ref{lem:structure_block}]
Clearly, $A$ has at least $2$ vertices.
If $A$ has $2$ vertices, then it is $K_2$ possibly with a unique single pendant edge attached.
If $A$ has at least $3$ vertices, then by Lemma~\ref{lem:essence}, $\overline{A}$ is $3$-connected or a cycle.
\end{proof}

\begin{proof}[Proof of Lemma~\ref{lem:structure_proper}]
Clearly $A$ has at least $3$ vertices.
By Lemma~\ref{lem:essence}, $\overline{A^+}$ is $3$-connected or a cycle.
\end{proof}

\begin{proof}[Proof of Lemma~\ref{lem:interiors}]
We have
$$A_1\cap A_2=(\int A_1\cap \int A_2)\cup (\int A_1\cap \partial A_2)\cup (\partial A_1\cap \int A_2)\cup (\partial A_1\cap \partial A_2).$$

\begin{mclaim}{A}
If $\int A_1\cap \int A_2\neq\emptyset$, then $\partial A_1\cap \int A_2\neq\emptyset$ and $\partial A_1\cap \partial A_2\neq \emptyset$.
\end{mclaim}

\begin{proofclaim}
If a vertex $v\in \int A_1\cap \int A_2$, then $A_1$ ($A_2$) is a block atom, or a proper atom.
By connectedness there is a path $P$ in $A_1$ joining $v$ to $\partial A_1$.
Either all the internal vertices of $P$ belong to $\int A_1\cap \int A_2$, or one of them is in $\partial A_1\cap \int A_2$.
In the latter case, we are done.
In the former case, we have an edge $e\in \int A_1\cap \int A_2$ that is moreover incident to $\partial A_1$.
If $\int A_2\subseteq \int A_1$, then $A_2\subseteq A_1$.
By minimality, $A_2=A_1$, which is a contradiction. Hence there exists an edge $f\in \int A_2\setminus \int A_1$ incident to $e$ and both $e$ and $f$ are incident to the same vertex in $\partial A_1$.
It follows that $\partial A_1\cap\int A_2\neq\emptyset$.
By symmetry, we have $\partial A_2\cap\int A_1\neq\emptyset$, as well.
\end{proofclaim}

\begin{mclaim}{B}
If $\partial A_1\cap \int A_2\neq\emptyset$, or $\partial A_2\cap \int A_1\neq\emptyset$, then both $A_1$ and $A_2$ are block or proper atoms.
\end{mclaim}

\begin{proofclaim}
Assume $\partial A_1\cap \int A_2\neq\emptyset$.
Clearly, $A_2$ cannot be a star, or a dipole, since
these types of atom contain no vertices in the interior.
Assume $A_1$ is a star, or a dipole with $|\partial A_1\cap \int A_2|=1$.
Since $\partial A_1$ is a non-trivial articulation, $A_2$ is neither a block atom, nor a proper atom, thus we get a contradiction.
Let $A_1$ be a dipole and $|\partial A_1\cap \int A_2|=2$.
Since $\int A_2$ is a component of $X\setminus\partial A_2$, $\int A_2$ contains with $\partial A_1$ all the edges of $A_1$.
Hence $A_1\subseteq A_2$, which implies $A_1=A_2$, a contradiction.
By symmetry we get the same conclusion assuming $\partial A_2\cap \int A_1\neq\emptyset$.
\end{proofclaim}

\begin{mclaim}{C}
We have $\partial A_1\cap \int A_2=\emptyset$ and $\partial A_2\cap \int A_1=\emptyset$.
\end{mclaim}

\begin{proofclaim}
Assume to the contrary that $\partial A_1\cap \int A_2\neq\emptyset$.
By Claim~B, each of $A_1$ and $A_2$ is a block, or a proper atom.
We distinguish three cases.
\begin{packed_itemize}
\item $\int A_1 \cap \bo A_2 = \emptyset$.
Since $\int A_2$ is a connected component $X\setminus \bo A_2$, it contains together with $\partial A_1\cap \int A_2$ all vertices and edges of $\int A_1$.
We have $\int A_1 \subseteq \int A_2$, which is a contradiction.
\item $|\int A_1\cap \bo A_2| = 1$.
By Lemma~\ref{lem:structure_block} and~\ref{lem:structure_proper} we have that $A_1$ is $2$-connected and using similar argument as in the previous case we get that $\int A_1 \subseteq \int A_2$, which is a contradiction.
\item $|\int A_1\cap \bo A_2| = 2$.
We may assume that $\bo A_2$ is a $2$-cut in $\int A_1$.
It follows that $A_2$ is a proper atom.
By Lemma~\ref{lem:structure_block} and~\ref{lem:structure_proper}, we have that $A_1$ is either a block atom isomorphic to a cycle, or $A_1$ is a proper atom.
In the former case, both vertices of $\bo A_2$ are of degree $2$, which is a contradiction.
By the same argument we exclude the case when $A_1$ is a proper atom isomorphic to a path.

Since $\int A_2$ is a component of $X\setminus \bo A_2$, the two vertices of $\bo A_1$ are connected in $A_2$.
Since $\bo A_2$ is a $2$-cut in $X$, there is a component $B \neq \int A_2$ in $X\setminus \bo A_2$.
By the assumptions, $\bo A_1 \not\subseteq B$.
Hence, $B \subseteq \int A_1$ is a part, which is a contradiction.
\end{packed_itemize}
By symmetry we also derive a contradiction assuming $\partial A_2\cap \int A_1=\emptyset$.
\end{proofclaim}

We summarize the above claims as follows.
By Claim~C we have that  $\partial A_2\cap \int A_1=\emptyset$ and $\partial A_1\cap \int A_2=\emptyset$.
By Claim~A and Claim~C we have $\int A_1\cap \int A_2=\emptyset$, as well.
This completes the proof of the lemma.
\end{proof}

\begin{proof}[Proof of Lemma~\ref{lem:atoms_and_automorphisms}]
Every automorphism permutes the set of articulations and non-trivial $2$-cuts.
Thus, $g(\bo A)$ separates $g(\int A)$.
\end{proof}

\section{Babai's theorem}
\label{app:babai}

Since the Babai's paper~\cite{babai1975automorphism} might be difficult to find online, we decided to rewrite the full statement of his theorem which gives a description of $\Aut(\planar)$.

\begin{theorem}[Babai~\cite{babai1975automorphism}, 8.12 The Main Corollary] \label{thm:babai}
Let $\Psi$ be a finite group. All graphs below are assumed to be finite.
\begin{enumerate}[(A'')]
\item[(A)] $\Psi$ is representable by a planar graph if and only if
\begin{equation} \label{eq:babai_A}
\Psi \cong \Psi_1 \wr \gS_{n_1} \times \cdots \times \Psi_t \wr \gS_{n_t}
\end{equation}
for some $t$, $n_1,\dots,n_t$ where the groups $\Psi_i$ are representable by connected planar
graphs.
\item[(A')] $\Psi$ is representable by a planar graph with a fixed point if and only if $\Psi$ is
representable by a planar graph.
\item[(A'')] $\Psi$ is representable by a planar fixed-point free graph if and only if a
(\ref{eq:babai_A}) decomposition exists with all groups $\Psi_i$ possessing planar connected
fixed-point-free graph representation.
\item[(B)] For $|\Psi| \ge 3$, $\Psi$ is representable by a connected planar graph if and only if
\begin{equation} \label{eq:babai_B}
\Psi \cong (\Psi_1 \wr \gS_k) \wr (\Psi_2 |_\calK)
\end{equation}
for some positive integer $k$, where $\Psi_1$ should be representable by a connected graph having a
fixed point; and either $\Psi_2$ is representable by a 2-connected planar graph $G_2$, or $k \ge 2$
and $|\Psi_2| = 1$. In the former case, $\Psi_2 |_\calK$ denotes the not necessarily effective
permutation group, acting on some orbit $\calK$ of $\Aut(G_2)$.
\item[(B')] $\Psi$ is representable by a connected planar graph having a fixed point if and only if
a (\ref{eq:babai_B}) decomposition exists as described under (B) with either $|\Psi_2| = 1$ or $G_2$
a 2-connected planar graph with a fixed point.
\item[(B'')] $\Psi$ is representable by a connected fixed-point-free planar graph if and only if a
(\ref{eq:babai_B}) decomposition exists with $|\Psi_2| \ne 1$, $G_2$ fixed-point free (hence
$|\calK| \ge 2$).
\item[(C)] $\Psi$ is representable by a 2-connected planar graph if and only if
\begin{equation} \label{eq:babai_C}
\Psi \cong (\Psi_1 \wr \gS_k) \wr (\Psi_2 |_\calK)
\end{equation}
where $|\Psi_1| \le 2$; $\Psi_2$ is representable by some 2-connected planar graph $G_2$. If
$|\Psi_1| = k = 1$, $G_2$ should be 3-connected. $\Psi_2 |_\calK$ denotes the action of $\Psi_2
\cong \Aut(G_2)$, as a not necessarily effective permutation group, on $\calK$, an orbit of either
an ordered pair $(a,b)$ or of an unordered pair $\{a,b\}$ of adjacent vertices $a,b \in V(G_2)$.
\item[(C')] $\Psi$ is representable by a 2-connected planar graph with a fixed point or with an
invariant edge if and only if a (\ref{eq:babai_C}) decomposition exists such that $G_2$ has a fixed
point or an invariant edge.
\item[(C'')] $\Psi$ is representable by a 2-connected planar fixed-point-free graph if and only if a
(\ref{eq:babai_C}) decomposition exists with $G_2$ fixed-point-free.
\item[(D)] $\Psi$ is representable by a 3-connected planar graph if and only if $\Psi$ is isomorphic
to one of the finite symmetry groups of the 3-space:
\begin{equation} \label{eq:babai_D}
\begin{split}
&\gC_n,\quad \gD_n,\quad \gA_4,\quad \gS_4,\quad \gA_5,\\
&\gC_n \times \gC_2,\quad \gD_n \times \gC_2,\quad \gA_4 \times \gC_2,\quad \gS_4
\times \gC_2,\quad \gA_5 \times \gC_2.
\end{split}
\end{equation}
\item[(D')] $\Psi$ is representable by a 3-connected planar graph with a fixed point if and only if
$\Psi$ is a cyclic or a dihedral group.
\item[(D'')] $\Psi$ is representable by a 3-connected planar fixed-point-free graph if and only if
$|\Psi| \ge 2$ and $\Psi$ is one of the groups listed under (\ref{eq:babai_D}).
\end{enumerate}
\end{theorem}

It reads in a nutshell as follows.  The automorphism group of a $k$-connected planar graph ($k\leq
2$) is constructed by combining automorphism groups of smaller $k$-connected planar graphs with
stabilizers of $k$-connected planar graphs and automorphism groups of $(k+1)$-connected graphs.

The part (A) corresponds to Theorem~\ref{thm:aut_disconnected}.  The automorphism groups listed
in the part (D) are the spherical groups described in Section~\ref{sec:3connected} and are based on
the classical results from geometry.  Therefore, the novel parts are (B) and (C). 
It is not clear which groups are $\Psi_2 |_\calK,$ used in (\ref{eq:babai_B}) and (\ref{eq:babai_C}).

In principle, it would be possible to derive Theorem~\ref{thm:babai} from the described Jordan-like
characterization of Theorems~\ref{thm:stabilizers} and~\ref{thm:composition}, and the reader
can work out further details. The opposite is not possible because Jordan-like characterization
contains more information about automorphism groups of planar graphs.
\end{document}